\documentclass[letterpaper, 10 pt, conference]{ieeeconf}  

\IEEEoverridecommandlockouts                              
\usepackage{graphics} 
\usepackage{epsfig} 
\usepackage{times} 
\usepackage{amsmath} 
\usepackage{amssymb}  
\usepackage{xspace}
\usepackage{caption}
\usepackage{subcaption}
\usepackage{color}

\newcommand{\htwo}{H\textsubscript{2}\xspace}
\title{\LARGE \bf
Dynamic Optimization and Optimal Control of Hydrogen Blending Operations in Natural Gas Networks
}

\author{Saif R. Kazi$^1$, Kaarthik Sundar$^2$ and Anatoly Zlotnik$^1$
\thanks{*This study was funded by the U.S. Department of Energy's Advanced Grid Modeling (AGM) projects ``Joint Power System and Natural Gas Pipeline Optimal Expansion'' and ``Dynamical Modeling, Estimation, and Optimal Control of Electrical Grid-Natural Gas Transmission Systems'', as well as LANL Laboratory Directed R\&D Project ``Efficient Multi-scale Modeling of Clean Hydrogen Blending in Large Natural Gas Pipelines to Reduce Carbon Emissions''.  Research conducted at Los Alamos National Laboratory is done under the auspices of the National Nuclear Security Administration of the U.S. Department of Energy under Contract No. 89233218CNA000001.}
\thanks{$^1$Saif R. Kazi and Anatoly Zlotnik are in the Applied Mathematics \& Plasma Physics Group, Los Alamos National Laboratory, Los Alamos, NM 87545, USA
        {\tt\small \{skazi,azlotnik\}@lanl.gov}}%
\thanks{$^2$Kaarthik Sundar is in the Information Systems \& Modeling Group, Los Alamos National Laboratory, Los Alamos, NM 87545, USA
        {\tt\small kaarthik@lanl.gov}}%
}

\begin{document}

\maketitle
\thispagestyle{empty}
\pagestyle{empty}

\begin{abstract}

We present a dynamic model for the optimal control problem (OCP) of hydrogen blending into natural gas pipeline networks subject to inequality constraints. The dynamic model is derived using the first principles partial differential equations (PDEs) for the transport of heterogeneous gas mixtures through long distance pipes. Hydrogen concentration is tracked together with the pressure and mass flow dynamics  within the pipelines, as well as mixing and compatibility conditions at nodes, actuation by compressors, and injection of hydrogen or natural gas into the system or withdrawal of the mixture from the network. We implement a lumped parameter  approximation to reduce the full PDE model to a differential algebraic equation (DAE) system that can be easily discretized and solved using nonlinear optimization or programming (NLP) solvers. We examine a temporal discretization that is advantageous for time-periodic boundary conditions, parameters, and inequality constraint bound values.  The method is applied to solve case studies for a single pipe and a multi-pipe network with time-varying parameters in order to explore how mixing of heterogeneous gases affects pipeline transient optimization.


\end{abstract}

\section{Introduction} \label{sec:intro}


The transition away from reliance on fossil fuels towards renewable and sustainable energy systems has inspired the development of new technologies for the production, transport, and utilization of energy.  At present, a large fraction of primary energy consumption for heating and power is still sustained by the extraction, pipeline transport, and combustion of natural gas.  Because hydrogen (\htwo) has high heating value per mass, can be produced using electrolysis powered by renewable electricity, and results in no emissions other than water when burned, efforts are underway to blend hydrogen into existing natural gas pipelines so these systems can also support the energy transition.  


The physical and chemical differences between hydrogen and natural gas significantly affect pipeline flow transients \cite{chaczykowski2018gas,baker2023transitions}, energy capacity \cite{tabkhi2008mathematical,sodwatana2023optimization}, and economics \cite{haeseldonckx2007use,zlotnik2023effects}. In blending, time-varying hydrogen injections are made into a pipeline that mainly carries natural gas to consumers with variable consumption \cite{guandalini2017dynamic}.  The gas composition must be monitored and optimized to account for energy content and pipeline operating efficiency, and optimal control is needed to handle the amplified transients caused by complex flow physics \cite{Kai_and_Saif}.  There is a need to examine how injection of much lighter hydrogen into natural gas pipelines can be accounted for in dynamic optimization, both conceptually and computationally.


Previous academic and industrial research found challenges for hydrogen injection into natural gas pipelines \cite{melaina2013blending}, for both pipes \cite{miao2021long} and compressor turbomachinery \cite{schuster2020centrifugal}.  Recent studies focus on simulating the complex dynamics that arise \cite{zhang2022modelling,witkowski2018analysis}.  The development of optimal control methods for managing the mixing of heterogeneous physical flows on networks defined by nonlinear PDEs has received less attention, particularly in the case of inequality constraints.


In this study, we extend optimal control methods created for gas pipeline systems \cite{rachford2000optimizing,hari2021operation} to account for injection of hydrogen and natural gas (\htwo-NG) throughout a network with complex time-dependent boundary conditions (BCs) and inequality constraints. In addition to computational challenges of increasing complexity and numerical ill-conditioning, we address the conceptual challenge of a well-posed formulation in which consumptions are variable.  To fill the conceptual gap, we free certain BCs for flow, determine conditions on concentration compatibility, and maximize an objective that represents the economic value created by transporting the mixture from suppliers of \htwo or NG to consumers of energy.  To address the computational challenge, we employ methods for approximating the partial differential equation (PDE) system for gas pipeline flow with a differential algebraic equation (DAE) system, and then discretizing the latter to a nonlinear program (NLP).  The optimal control problem (OCP) is applied to manage simple time-varying boundary concentration changes in a single pipe subject to constant inequalities on state variables, resulting in highly non-trivial solutions.  We demonstrate generalizability using a test pipeline network with three compressor stations and a loop.  The optimal control formulation and modeling approach employed here have not been synthesized and applied previously to the authors' knowledge.


The rest of the paper is as follows.  In Section \ref{sec:dynamics}, we describe the physical modeling of natural gas and hydrogen mixtures flowing through a pipe, as well as handling of BCs, approximation by lumped parameters, and time discretization.  In Section \ref{sec:network} the modeling is extended to mixing throughout a pipeline network, including discussion of spatial discretization of the graph, nodal compatibility conditions, and energy content.  The optimal control problem is defined in Section \ref{sec:optimal} by adding an objective function that accounts for the value created by the flow allocation and the cost of operations.  Section \ref{sec:results} provides implementation details and results of case studies for a single pipe and a test network, and we conclude in Section \ref{sec:conclusions}.

\section{Dynamic Heterogeneous Flow in a Pipe}  \label{sec:dynamics}


Isothermal transport of an \htwo-NG mixture in a pipe without transients that create waves or shocks can be approximated by the system of partial differential equations
\begin{subequations}\label{eq:full-PDE}
\begin{align}\label{eq:h2_continuity}
   & \frac{\partial \rho_{H_2}}{\partial t} +   \frac{\partial \phi_{H_2}}{\partial x} = 0, \\
\label{eq:ng_continuity}
   & \frac{\partial \rho_{NG}}{\partial t} +   \frac{\partial \phi_{NG}}{\partial x} = 0, \\
\label{eq:momentum}
   & \frac{\partial \phi}{\partial t} + \frac{\partial P}{\partial x} = -\frac{\lambda}{2 D} \frac{\phi |\phi|}{\rho},
\end{align}
\end{subequations}
where $\rho_{H_2}$ and $\rho_{NG}$ denote the partial densities, and $\phi_{H_2}$ and $\phi_{NG}$ denote the partial mass flux, of \htwo and NG respectively \cite{chaczykowski2018gas}. Here $P$ denotes the total pressure of the mixture and $\lambda$ and $D$ are friction factor and pipe diameter parameters. 

\subsection{Assumptions and Simplification} \label{sec:assume}

The total density $\rho $ and flux $\phi$ are assumed to be sums of component-wise partial densities $\rho_k$ and fluxes $\phi_k$:
\vspace{-1ex}
\begin{subequations}
\begin{align}
   & \rho_{H_2}  + \rho_{NG}  = \rho, \\
   & \phi_{H_2} + \phi_{NG} = \phi.
\end{align}
This holds true for ideal gas mixtures with additive partial pressures, and enables simplified PDE system modeling.
The ideal gas assumption relates the total pressure in terms of partial densities as:
\vspace{-1ex}
\begin{align}
    & P = p_{H_2} + p_{NG}, \\
    & P = a_{H_2}^2 \rho_{H_2} +  a_{NG}^2 \rho_{NG},
\end{align}
where $a_{H_2}$ and $a_{NG}$ are speed(s) of sound in pure \htwo and NG respectively. For ideal gases, the speed of sound is given by $a=\sqrt{\frac{RT}{M}}$ where $R$, $T$, and $M$ are universal gas constant, temperature, and molecular mass.  The equation of state is
\begin{equation} \label{eq:eos1}
    P = \left(a_{H_2}^2 \frac{\rho_{H_2}}{\rho_{H_2} + \rho_{NG}} +  a_{NG}^2 \frac{\rho_{NG}}{\rho_{H_2} + \rho_{NG}}\right) (\rho_{H_2} + \rho_{NG}). 
\end{equation}
The ratio of partial density to the total density, $\rho_{H_2}/(\rho_{H_2} + \rho_{NG})$, is the mass fraction of \htwo in the gas mixture, which we denote by $\eta$, and write equation \eqref{eq:eos1} as
\begin{equation}\label{eq:EOS2}
    P = \left(a_{H_2}^2 \eta +  a_{NG}^2 (1-\eta)\right) \rho = a^2 \rho.
\end{equation}
\end{subequations}
We make another assumption regarding the flux derivative $\partial \phi/\partial t$ to provide for numerical stability. In most practical cases, this term is negligible compared to the pressure gradient term $\partial P/\partial x$ \cite{osiadacz1984simulation}. For ideal gas mixtures, the ratio of these two terms is approximately $\mathcal{O}(u/a)$ where $u/a$ is the ratio of gas velocity and sound speed in the mixture. Typical gas velocity values of $u \in [1,10]$ m/s are much smaller than the speed of sound $a \in [ 350, 1000]$ m/s.
Using the above assumptions, we simplify the PDE system \eqref{eq:full-PDE} to the system
\vspace{-1ex}
\begin{subequations}\label{eq:simple-PDE}
\begin{align}\label{eq:h2_continuity_simp}
   & \frac{\partial \rho_{H_2}}{\partial t} +  \frac{\partial}{\partial x} (\eta \phi) = 0, \\
\label{eq:ng_continuity_simp}
   & \frac{\partial \rho_{NG}}{\partial t} +  \frac{\partial}{\partial x} ((1-\eta) \phi) = 0, \\
\label{eq:momentum_simp}
   &  \frac{d a^2 \rho}{d x} = -\frac{\lambda}{2 D} \frac{\phi |\phi|}{\rho}.
\end{align}
\end{subequations}

\subsection{Lumped Element Model and Time Discretization} \label{sec:lump}

For a short pipe of length $L$, the system can be integrated over the length variable from $x=0$ to $x=L$ to yield:
\begin{subequations}\label{eq:single-pipe-ODE}
\begin{align}\label{eq:h2_continuity_lumped}
   & L \frac{d \bar \rho_{H_2}}{d t} + \eta^L \phi^L - \eta^0 \phi^0  = 0, \\
\label{eq:ng_continuity_lumped}
   & L \frac{d \bar \rho_{NG}}{d t} + (1-\eta^L) \phi^L - (1-\eta^0) \phi^0  = 0, \\
\label{eq:momentum_lumped}
   & a^2(L) \rho (L) - a^2 (0) \rho (0)  = -\frac{\lambda L}{2 D} \frac{\bar \phi |\bar \phi|}{\bar \rho}.
\end{align}
The partial densities and the sound speed are evaluated at nodes whereas the flux and concentration variables are defined at the endpoints of the pipe edge. The right side of equation \eqref{eq:momentum} is approximated using the variables $\bar \rho = (\rho(L) + \rho(0))/2$ and $\bar \phi = (\phi^L + \phi^0)/2$, which denote the average values for density and flux respectively. 



\end{subequations}

A simple first order forward finite difference formula is used for time discretization on the interval $[0,T^f]$ using $N$ points $T^S=\{t_n\}_{n=1}^N$ defined by $t_n= \Delta t(n-1)$ for $n=1,2,\ldots,N$.  The resulting approximation is
\begin{subequations}
\begin{align}
    \!\!\left.\frac{d \bar \rho_{H_2}}{d t}\right|_{t_n} & \!\!\!\!\! \approx \! \frac{\left(\!\left.\rho^L_{H_2}\right|_{t_{n+1}} \!\!\!\!\!\! + \left.\rho^0_{H_2}\right|_{t_{n+1}}\!\right) \!-\! 
    \left(\!\left.\rho^L_{H_2}\right|_{t_{n}} \!\!\!\! + \left.\rho^0_{H_2}\right|_{t_{n}}\!\right)}{2\Delta t}   \\
   \!\! \left.\frac{d \bar \rho_{NG}}{d t}\right|_{t_n} & \!\!\!\!\! \approx \! \frac{\left(\!\left.\rho^L_{NG}\right|_{t_{n+1}} \!\!\!\!\!\! + \left.\!\rho^0_{NG}\right|_{t_{n+1}}\!\right) \!-\! 
    \left(\!\left.\rho^L_{NG}\right|_{t_{n}} \!\!\!\! + \left.\!\rho^0_{NG}\right|_{t_{n}}\!\right)}{2\Delta t}  
\end{align}
\end{subequations}

\subsection{Boundary Conditions} \label{sec:bound}

Instead of specifying initial and terminal conditions in the OCP, we impose cyclic BCs on the state variables, of form
\begin{equation}\label{eq:cyclic_BC}
    x(T^f) - x(0) = 0,
\end{equation}
where $x$ denotes any of the partial densities $\rho_{H_2}$ or 
$\rho_{NG}$ and $T^f$ is the final time point. Rather than explicitly enforcing the boundary condition equation \eqref{eq:cyclic_BC}, we rewrite the time derivative for the final time step as 
\begin{equation}
    \left.\frac{d \bar x}{dt}\right|_{t_{N}}  \approx \frac{x(0) - x(t_N)}{\Delta t}.
\end{equation}
This formulation implicitly includes the cyclic condition \eqref{eq:cyclic_BC} and reduces the number of variables and constraints, which improves the numerical conditioning of the resulting NLP.  Moreover, such elliptic BCs result in very well-posed problems that yield well-behaved solutions.  The method can be adapted to non-periodic BCs by extending the time horizon and interpolating, which can be applied to real data in a model-predictive moving horizon manner, as was previously demonstrated for transport of a homogeneous gas \cite{rudkevich2019evaluating}.

\subsection{Non-dimensional system of equations}  \label{sec:nondim}

As is common practice \cite{herty2010new}, the variables and equations in the model are non-dimensionalized for numerical purposes. Standard or nominal values for length, pressure, and velocity are chosen and denoted by $l_0, p_0 \text{ and } \eta_0 = 1$, respectively. Nominal values are also derived for wave speed, $a_0 = \sqrt{a_{H_2} \cdot a_{NG}} \approx 672 \text{ m s}^{-1}$, and flow speed, $v_0 = a_0/M$, where $M$ denotes the Mach value for gas velocity. For the purpose of this study, we use a nominal value of $M = 1/300$, which results in $v_0 \approx 2.24 \text{ m s}^{-1}$.  Nominal density is chosen as $\rho_0 = p_0/a^2_0$ and nominal flow and flux are $f_0  = \rho_0 A_0 v_0$ and $\phi_0 = \rho_0 v_0$, where nominal pipe cross section area is $A_0 = 1$. Setting $\hat \rho = \rho/\rho_0, \hat a = a/a_0, \hat L = L/l_0, \hat D = D/l_0 \text{ and } \hat \phi = \phi/\phi_0$, we find that equations (\ref{eq:single-pipe-ODE}) reduce to
\vspace{-2.5ex}
\begin{subequations}\label{eq:dimensionless-single-pipe-ODE}
\begin{align}\label{eq:h2_continuity_nd}
   & \hat L \frac{d \hat{\bar{\rho_{H_2}}}}{d t} + \frac{\eta_0 \phi_0}{\rho_0 l_0} \left(\hat \eta^{\hat L} \hat \phi^{\hat L} - \hat \eta^0 \hat \phi^0 \right) = 0, \\
\label{eq:ng_continuity_nd}
   & \hat L \frac{d \hat{\bar{\rho_{NG}}}}{d t} + \frac{\eta_0 \phi_0}{\rho_0 l_0} \left((1-\hat \eta^{\hat L}) \hat \phi^{\hat L} - (1-\hat \eta^0) \hat \phi^0 \right)  = 0, \\
\label{eq:momentum_nd}
   & \hat a^2(\hat L) \hat \rho (\hat L) - \hat a^2 (0) \hat \rho (0)  = - \left(\frac{\phi^2_0}{a^2_0 \rho^2_0}\right) \left(\frac{\lambda \hat L}{2 \hat D}\right) \frac{\hat{\bar{\phi}}|\hat{\bar{\phi|}}}{\hat{\bar{\rho}}}.
\end{align}
\end{subequations}
\vspace{-1ex}
By defining
\begin{equation}
    \kappa \triangleq \frac{\eta_0 \phi_0}{\rho_0 l_0} \!=\! \frac{v_0}{l_0} \,\, \text{and} \,\, \beta \triangleq \left(\frac{\phi^2_0}{a^2_0 \rho^2_0}\cdot\frac{\lambda \hat L}{2 \hat D}\right) \!=\! \frac{1}{M^2} \left(\frac{\lambda \hat L}{2 \hat D}\right),
\end{equation}
the system is rewritten in terms of the flow $\hat f = \hat \phi \hat A$ as
\begin{subequations}\label{eq:beta-dimensionless-single-pipe-ODE}
\begin{align}\label{eq:h2_continuity_ndrw}
   &\frac{d \hat{\bar{\rho_{H_2}}}}{d t} + \frac{\kappa}{\hat L \hat A} \left(\hat \eta^{\hat L} \hat f^{\hat L} - \hat \eta^0 \hat f^0 \right) = 0, \\
\label{eq:ng_continuity_ndrw}
& \frac{d \hat{\bar{\rho_{NG}}}}{d t} + \frac{\kappa}{\hat L \hat A} \left((1-\hat \eta^{\hat L}) \hat f^{\hat L} - (1-\hat \eta^0) \hat f^0 \right)  = 0, \\
\label{eq:momentum_ndrw}
  &   \hat a^2(\hat L) \hat \rho (\hat L) - \hat a^2 (0) \hat \rho (0)  = - \beta \frac{\hat{\bar{\phi}}|\hat{\bar{\phi|}}}{\hat{\bar{\rho}}}.
\end{align}
\end{subequations}

\vspace{-1ex}
\noindent In subsequent discussions, for ease of presentation, we omit the hat symbol that designates non-dimensional quantities with the understanding that all quantities are dimensionless.

\vspace{-0.5ex}
\section{Dynamic Heterogeneous Flow in a Network}  \label{sec:network}

The network model consists of pipes and compressors that connect nodes that are subject to injection and withdrawal flows, or non-optimized flows. The following graph-based notations are used to formulate the network-wide dynamic model for optimal transport of mixed gas. Set notations are:

\begin{itemize}
    \item $\mathbb G  = (\mathbb N, \mathbb P \cup \mathbb C)$ - graph of a pipeline network 
    \item $\mathbb N, \mathbb P, \mathbb C$ - sets of nodes, pipes, compressors, respectively
    \item $\mathbb S, \mathbb I, \mathbb W$ - sets of slack, injection and withdrawal nodes
    \item $(i,j) \in \mathbb P$ - the pipe connecting nodes $i$ and $j$
    \item $(i,j) \in \mathbb C$ - the compressor between nodes $i$ and $j$
\end{itemize}
\vspace{1ex}

\noindent Nodal variable notations are:
\begin{itemize}
    \item $\eta_i$ - concentration of \htwo (after mixing) at node $i$
    \item $a_i$ - speed of sound in gas at node $i$
    \item $\rho_{H_2,i},\rho_{NG,i}$ - partial density of \htwo or NG at node $i$
    \item $q^s_i,q^w_i$ - supply and withdrawal flows at injection $i \in \mathbb I$ and withdrawal $i \in \mathbb W$ nodes, respectively
    \item $g^E_i$ - energy delivered to withdrawal node $i \in \mathbb W$
\end{itemize}
\noindent Edge variable notations are:
\begin{itemize}
    \item $f^0_{ij},f^{L}_{ij}$ - flow through pipe $(i,j)$ at end points
    \item $\gamma^0_{ij},\gamma^L_{ij}$ - concentration of \htwo in pipe $(i,j)$ at end points
    \item $f^c_{ij}$ - flow through a compressor $(i,j)$
    \item $\alpha_{ij}$ - compressor ratio in compressor $(i,j)$
\end{itemize}

\noindent Parameter notations are:
\begin{itemize}
    \item $\lambda_{ij}$ - friction factor of pipe $(i,j)$
    \item $A_{ij}$ - cross-sectional area of pipe $(i,j)$
    \item $L_{ij,}D_{ij}$ - length and diameter of pipe $(i,j)$
    \item $\beta_{ij}$ - effective resistance of pipe $(i,j)$
    \item $\eta^s_i$ - concentration of injection at node $i \in \mathbb I \cup \mathbb S$
\end{itemize}



\noindent In addition to the dynamics on each pipe in equation \eqref{eq:beta-dimensionless-single-pipe-ODE}, we include algebraic equalities and inequality constraints in the developed OCP.  We suppose that a compressor station $(i,j)\in\mathbb{C}$ has pressure boost ratio $\alpha_{ij}$ and flow constraints
\vspace{-.5ex}
\begin{subequations} \label{eq:comp_constr}
\begin{align} 
    &p_{j}^{2} = \alpha_{ij}^{2} p_{i}^{2},  & \forall (i,j) \in \mathbb{C}, \label{eq:comp_boost} \\
    & 1\leq \alpha_{ij} \leq \alpha_{ij}^{\max}, & \forall (i,j) \in \mathbb{C}, \label{eq:comp_limits} \\
    &0 \leq f_{ij}^c \leq f_{ij}^{c,\max},    &\forall (i,j) \in \mathbb{C}. \label{eq:comp_flow} 
\end{align} 
\end{subequations}

\vspace{-.5ex}
\noindent We enforce flow balance for hydrogen and natural gas by
\begin{subequations}\label{eq:nodebalance}
\begin{align} \label{eq:h2flowbalance}
    & \!\!\!\!\!\! \sum_{i\in \partial_j^+} \gamma^L_{ij} f_{ij} - \sum_{k\in \partial_j^-} \gamma_{ij}^{0} f_{jk} =  \eta^s_{j} q^s_j -   \eta_{j} q^w_j, \,\, \forall j\in\mathbb{N}, \\ \label{eq:ngflowbalance}
    & \!\!\!\!\!\! \sum_{i\in \partial_j^+} (1-\gamma^L_{ij}) f_{ij} - \sum_{k\in \partial_j^-} (1-\gamma_{ij}^{0}) f_{jk} \nonumber \\ 
    & \!\!\!\!\!\! \qquad \qquad \qquad =  (1-\eta^s_{j}) q^s_j -  (1-\eta_{j}) q^w_j, \,\, \forall j\in\mathbb{N}.
\end{align} 
\end{subequations}

\vspace{-.5ex}
\noindent Here \(\partial_{j}^{+}\) and \(\partial_{j}^{-}\) are the sets of nodes connected to node $j$ by incoming and outgoing edges, respectively.  We suppose that only one of $q_j^s\geq 0$ and $q_j^w\geq 0$ can be positive at any given time for any $j$. Adding together equations \eqref{eq:h2flowbalance} and \eqref{eq:ngflowbalance} yields the total nodal mass balance.  
We define the energy content $g^E$ of the mixture withdrawal flow as
\begin{equation}\label{eq:energy_withdrawal}
    g^E_i = (\eta_i R_{H_2} + (1-\eta_i) R_{NG})q^w_i, \quad \forall i \in \mathbb W,
\end{equation}

\vspace{-.5ex}
\noindent where $R_{H_2}$ and $R_{NG}$ are calorific values for \htwo and NG, respectively.
Instead of specifying a fixed value or profile for withdrawal demand flow $q^w$, we either specify a fixed energy withdrawal $\bar g^E$ or impose an upper bound on the non-negative variable as maximum energy demand by
\begin{equation}\label{eq:energy_withdrawal_bound}
    0 \leq g^E_i \leq g^{E,max}_i, \quad \forall i \in \mathbb W.
\end{equation}

\vspace{-.5ex}
\noindent We also impose equations for all $t \in T^s$ that enforce compatibility conditions, determine the local sound speed, and delimit pressure within acceptable levels.  These are
\begin{subequations}\label{eq:network_equations}
\begin{equation}\label{eq:node-concentration}
    \eta_i = \frac{\rho_{H_2,i}}{(\rho_{H_2,i} + \rho_{NG,i})}, \quad \forall \, i \in \mathbb{N},
\end{equation}
\begin{equation}\label{eq:node-edge continuity}
    \eta_i - \gamma^0_{ij}= 0, \quad \forall \, (i,j) \in \mathbb{P}\cup\mathbb{C},
\end{equation}
\begin{equation} \label{eq:node-sound-speed}
    a^2_i = a^2_{H_2}\eta_{i} + a^2_{NG}(1-\eta_{i}), \quad \forall \, i \in \mathbb{N},
\end{equation}
\begin{equation} \label{eq:node-slack-pressure}
    a^2_i (\rho_{H_2,i} + \rho_{NG,i}) = p_i^s \,\, (\text{fixed}), \quad \forall \, i \in \mathbb{S},
\end{equation}
\begin{equation} \label{eq:node-pressure-bounds}
    p^{min}_i \leq a^2_i (\rho_{H_2,i} + \rho_{NG,i}) \leq p^{max}_i, \quad \forall \, i \in \mathbb{N}.
\end{equation}
\end{subequations}
Equations \eqref{eq:node-concentration}-\eqref{eq:node-pressure-bounds} define nodal concentration, node-edge concentration continuity at pipe inlet, nodal sound speed, slack node pressure, and nodal pressure bounds, respectively.

Empirical studies have shown \cite{osiadacz1984simulation,grundel2014model} that the lumped model approximation in equation \eqref{eq:single-pipe-ODE} is acceptable for short pipes ($< \approx$10km) in the usual transient pressure and flow regime seen in gas transmission pipelines.  Thus, we discretize pipelines in the network uniformly into short pipes of fixed length ($\Delta L$) as shown in Fig. \ref{sinlge-pipe} below.
\vspace{-2ex}
\begin{figure}[h!]
    \centering
    \includegraphics[height=20mm,width=80mm]{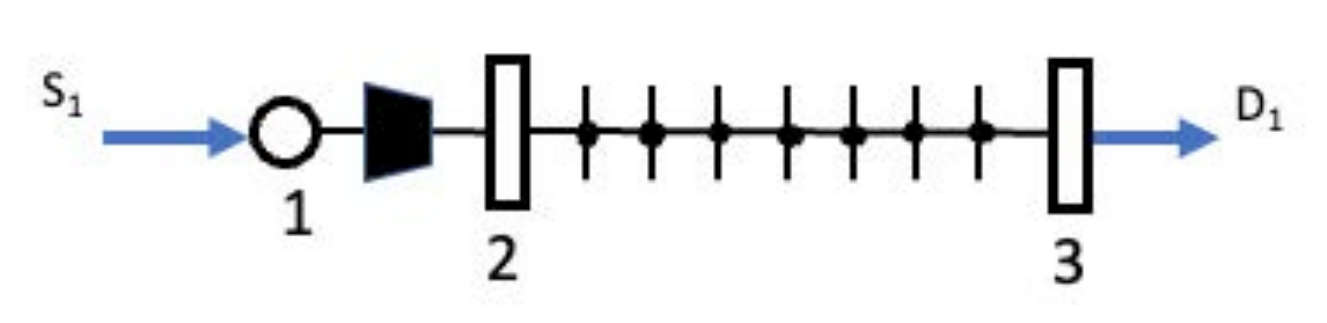}
    \vspace{-2ex}
    \caption{Single Discretized Pipe}
    \label{sinlge-pipe}
\end{figure}
\vspace{-3ex}

\section{Optimal Control Problem} \label{sec:optimal}


The objective function for the optimal control problem (OCP) consists of the economic value generated by the pipeline operator through sales of energy and purchases of gas constituents, minus the operating cost of compressors. The economic objective function is of the form
\begin{subequations}
\begin{equation}\label{eq:gas_sales_revenue}
\begin{split}
    R_e \!=\!\!  \sum_{t \in T_s} \left(\sum_{i \in I} (\eta^s_i c_{H_2,i} \!+\! (1-\eta^s_i) c_{NG,i}) q_{i,t}^s
    \!-\! \sum_{i \in W} C^E g_{i,t}^E \right)\!\!,
\end{split}
\end{equation}
where $c_{H_2}$ and $c_{NG}$ denote supplier offer prices for \htwo and NG, and $C^E$ denotes the price that consumers bid for energy delivered through the mixed gas in \$/MW.
The operating cost for each compressor is determined using an adiabatic compression work expression \cite{Sundar2018}, of form
\begin{equation}\label{eq:adiabatic_compression}
    W_c = \left( \frac{286.76\cdot \mu_{ij}\cdot T}{G_{ij}{(\mu_{ij}-1)}} \right) \left( \alpha_{ij}^{m_{ij}}-1 \right) \left|f^c_{ij}\right|,
\end{equation}

\vspace{-1ex}
\noindent
where $m_{ij}=(\mu_{ij}-1)/\mu_{ij}$ with $\mu_{ij}$ and $G_{ij}$ denoting the specific heat capacity ratio and specific gravity of the blend that may vary with location and time.  The above equation \eqref{eq:adiabatic_compression} is a simplification of more detailed adiabatic compression modeling \cite{menon05}.  Because $\mu_{ij}$ and $G_{ij}$ do not vary significantly for the hydrogen fractions considered here, we approximate them as constants for simplicity and denote the entire expression by $K$. We use nominal values $\mu = 1.31$  and $G = 0.505$, which correspond to 10\% \htwo by mass, and $T=288.7$ K as the compressor suction temperature.  We further approximate $m=1/2$, so that the cost of operating compressors is obtained by adding the total compressor work given by equation \eqref{eq:adiabatic_compression} for each compressor and multiplying by a constant electricity price of $\zeta = 0.07$ \$/kWh, yielding
\begin{equation}\label{eq:compressor_work}
    R_c =  K \cdot \zeta \cdot \sum_{t_n} \sum_{ij \in C} f^{c,t}_{ij} \left(\sqrt{\alpha^t_{ij}}-1\right).
\end{equation}

\vspace{-1ex}
\noindent
We use $m=1/2$ in the objective function because the economic value $R_e$ in \eqref{eq:gas_sales_revenue}  of operating the pipeline is expected to be about two orders of magnitude greater than the cost $R_c$ in \eqref{eq:compressor_work} of operating the pipeline by gas compression, so that the latter quantity does not significantly affect the objective value but is nonetheless useful to promote well-posedness and regularization.  The total objective function (to be minimized) is a weighted sum of the economic cost and the operating cost with a weight $\xi$, of form
\begin{equation}\label{eq:total_objective}
    \min \enskip\mathbb{T^{\text{cost}}} = \xi R_e + (1-\xi) R_c.  
\end{equation}

\vspace{-1ex}
\noindent
The scaling parameter can be used to prioritize on maximizing gas delivery or minimizing the operating cost. For $\xi = 0.5$, the objective reduces to the sum of the two terms.
\end{subequations}


\noindent The complete optimization problem is formulated as:
\vspace{-0.5ex}
\begin{flalign*} 
    \mathtt{OPT} ~ \triangleq ~ & \min ~~ \mathbb{T^{\text{cost}}} ~ \text{in } \eqref{eq:total_objective}, ~ \text{ -- subject to: } \\
    & \text{equations \eqref{eq:beta-dimensionless-single-pipe-ODE}, \eqref{eq:comp_constr}} \text{ -- (transport \& compression)} \\
    & \text{equations \eqref{eq:nodebalance}} \text{ -- (\htwo and NG nodal balance)} \\
    & \text{equations \eqref{eq:energy_withdrawal}, \eqref{eq:energy_withdrawal_bound}} \text{ -- (energy calculation)} \\
    & \text{equations \eqref{eq:network_equations}} \text{ -- (network compatibility)}
\end{flalign*}

\vspace{-0.5ex}
\noindent
and the optimized variables are all nodal and edge variables as defined in Section \ref{sec:dynamics}.

\section{Results}  \label{sec:results}

We solve the above problem for two examples that consist of a single pipe and an 8-node test network subject to variable boundary conditions and/or constraint bounds. 

\subsection{Computational Implementation} \label{sec:implement}

Problem $\mathtt{OPT}$ is solved for both test cases using a NLP solver such as KNITRO \cite{knitro} in the Julia based modeling language JuMP \cite{jump}. The problems are solved for a 24-hour time horizon with $T_f = 24$ hours in a sequential manner. At first, the steady-state problem is solved by setting $\Delta t = T_f$ so that the time discretization is a single point. Thereafter, the steady-state solution is used as the initial guess for the full time discretization with $\Delta t = 0.5 \text{ or } 1.0$ hours.  Note that because the problem is non-convex, there are no guarantees of global optimality for the solution obtained by the NLP solver.  Our approach ensures a feasible solution when a local optimum is attained, which is challenging for this problem.

\subsection{Single Pipe Case Study}    \label{sec:singlepipe}

\begin{figure}[h!]
    \centering
    \begin{subfigure}[b]{1.0\linewidth}
        \centering
        \includegraphics[height=40mm,width=0.9\linewidth]{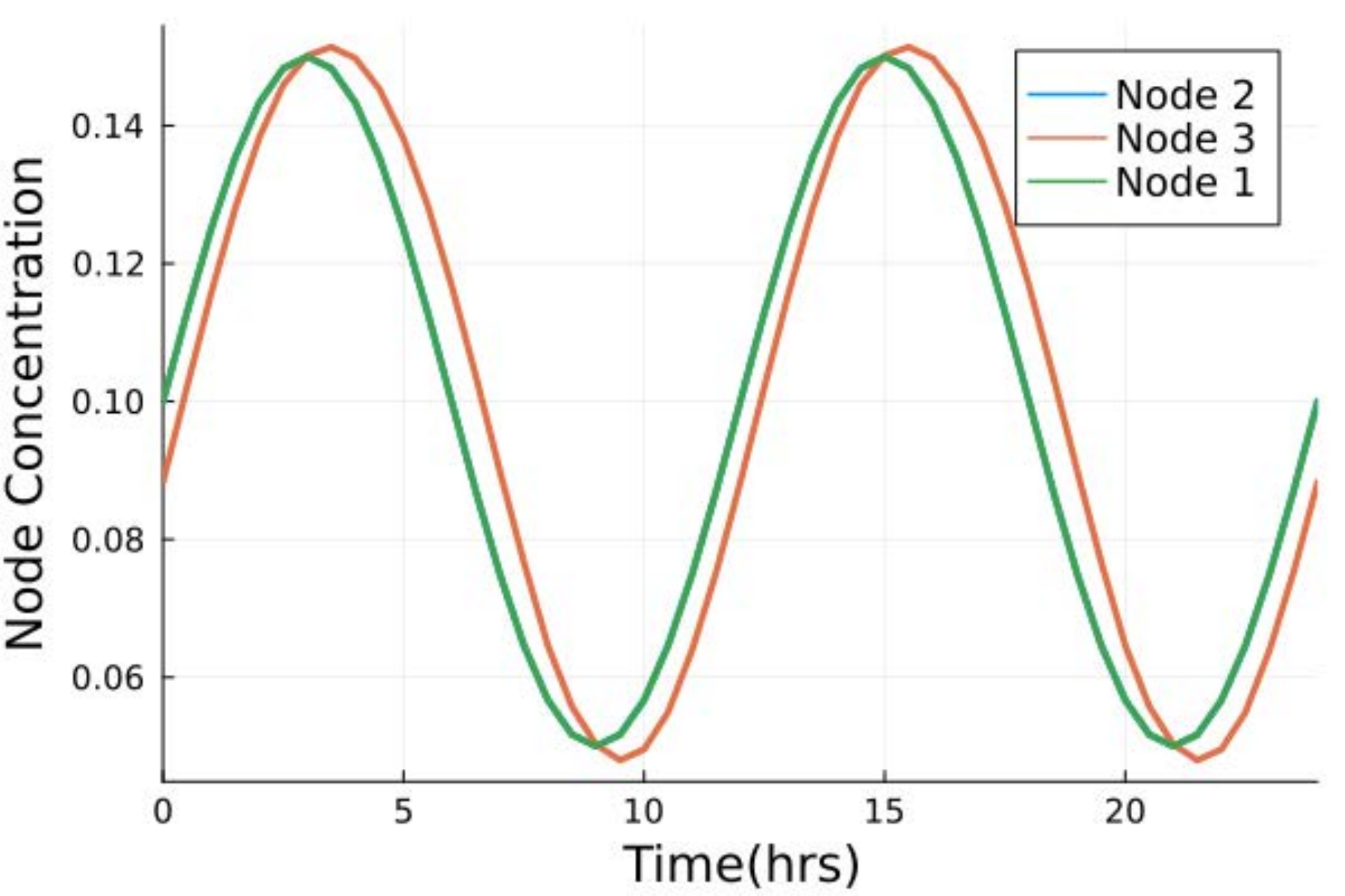}
        \caption{Nodal concentration}
        \label{Node concentration}
    \end{subfigure}
    \hfill
    \begin{subfigure}[b]{1.0\linewidth}
        \centering
        \includegraphics[height=40mm,width=0.9\linewidth]{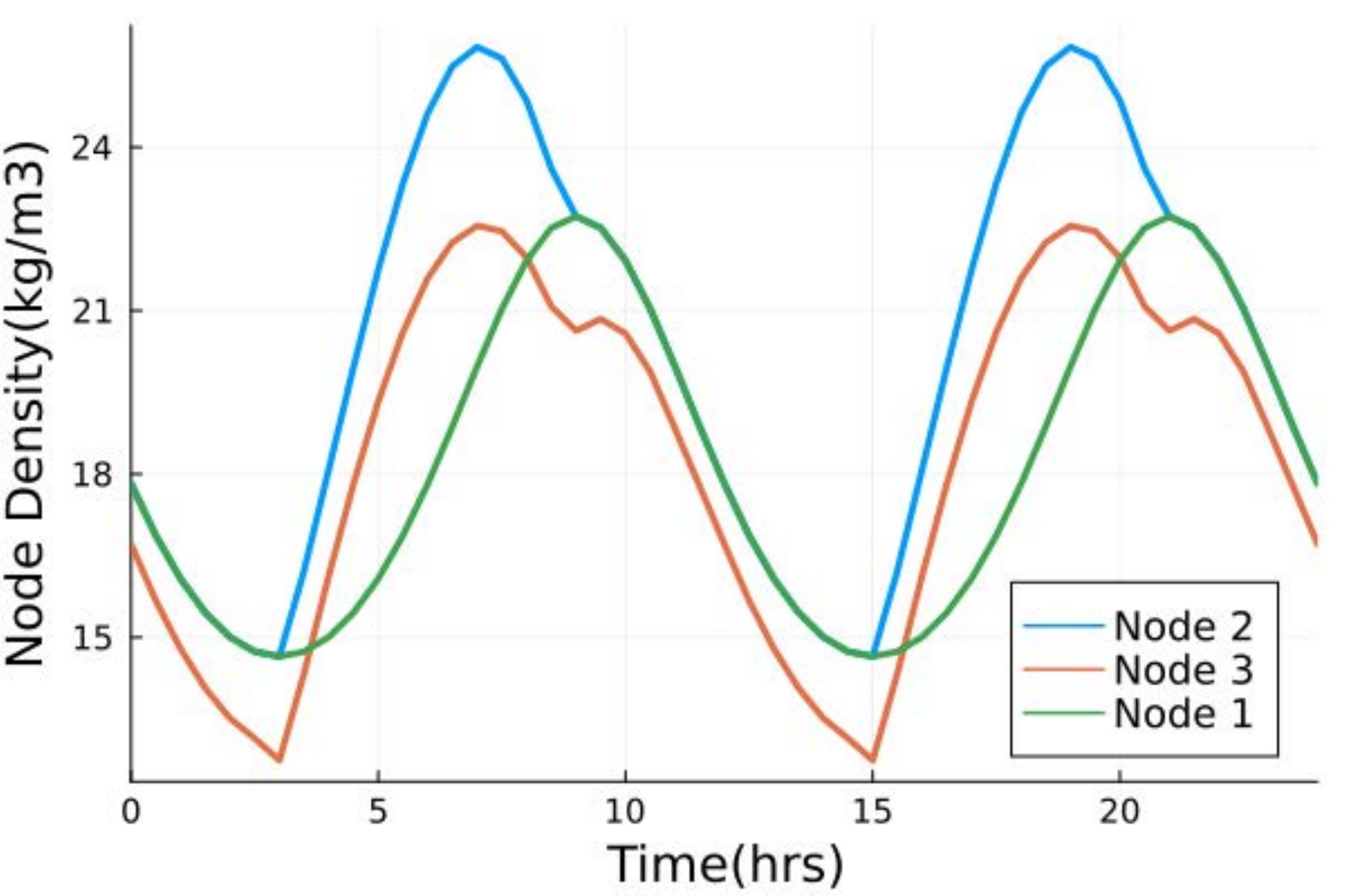}
        \caption{Nodal density}
        \label{Node density}
    \end{subfigure}
    \hfill
    \begin{subfigure}[b]{1.0\linewidth}
        \centering
        \includegraphics[height=40mm,width=0.9\linewidth]{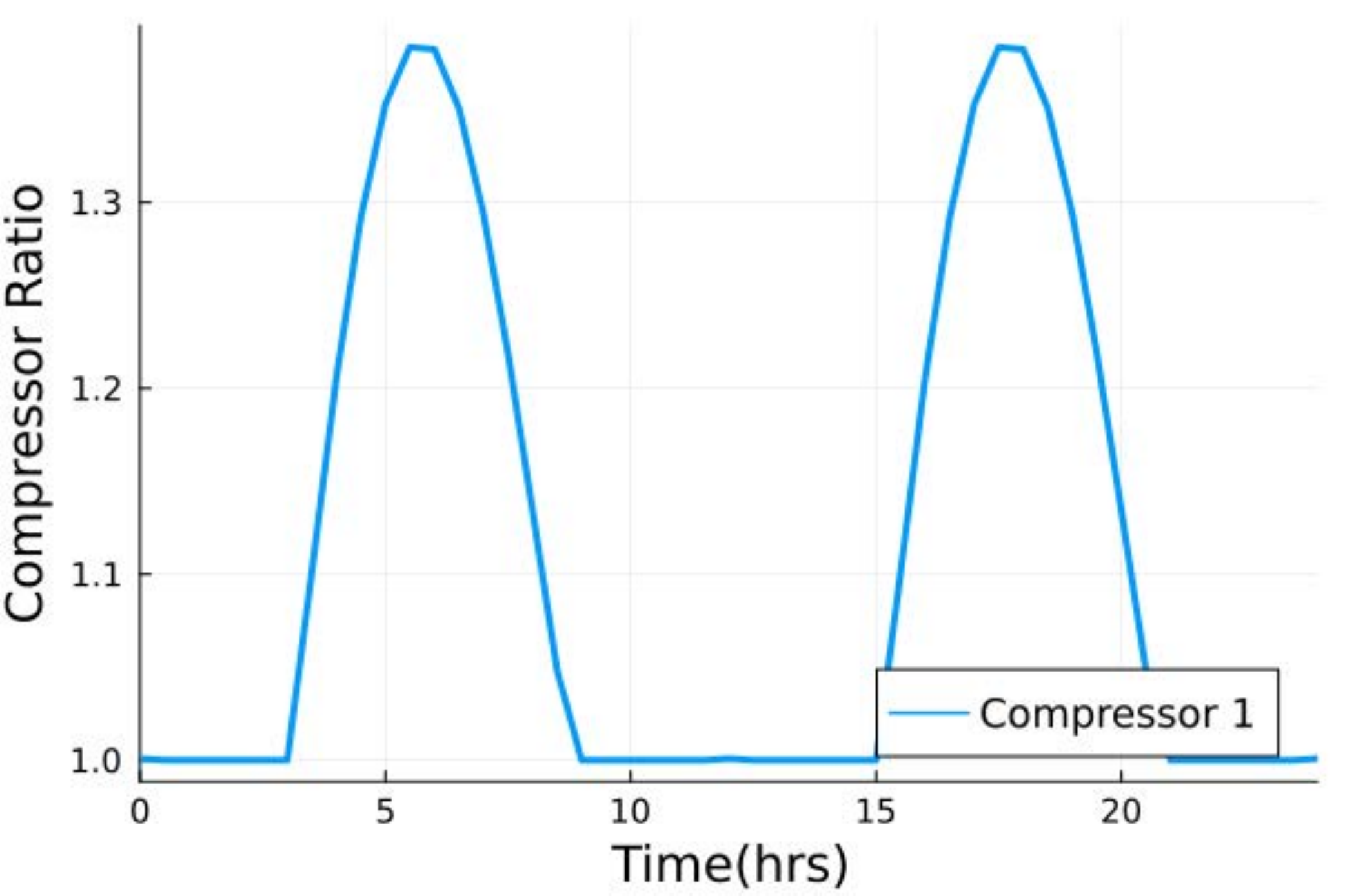}
        \caption{Compressor Ratio}
        \label{Compressor Ratio}
    \end{subfigure}
    \hfill
    \begin{subfigure}[b]{1.0\linewidth}
        \centering
        \includegraphics[height=40mm,width=0.9\linewidth]{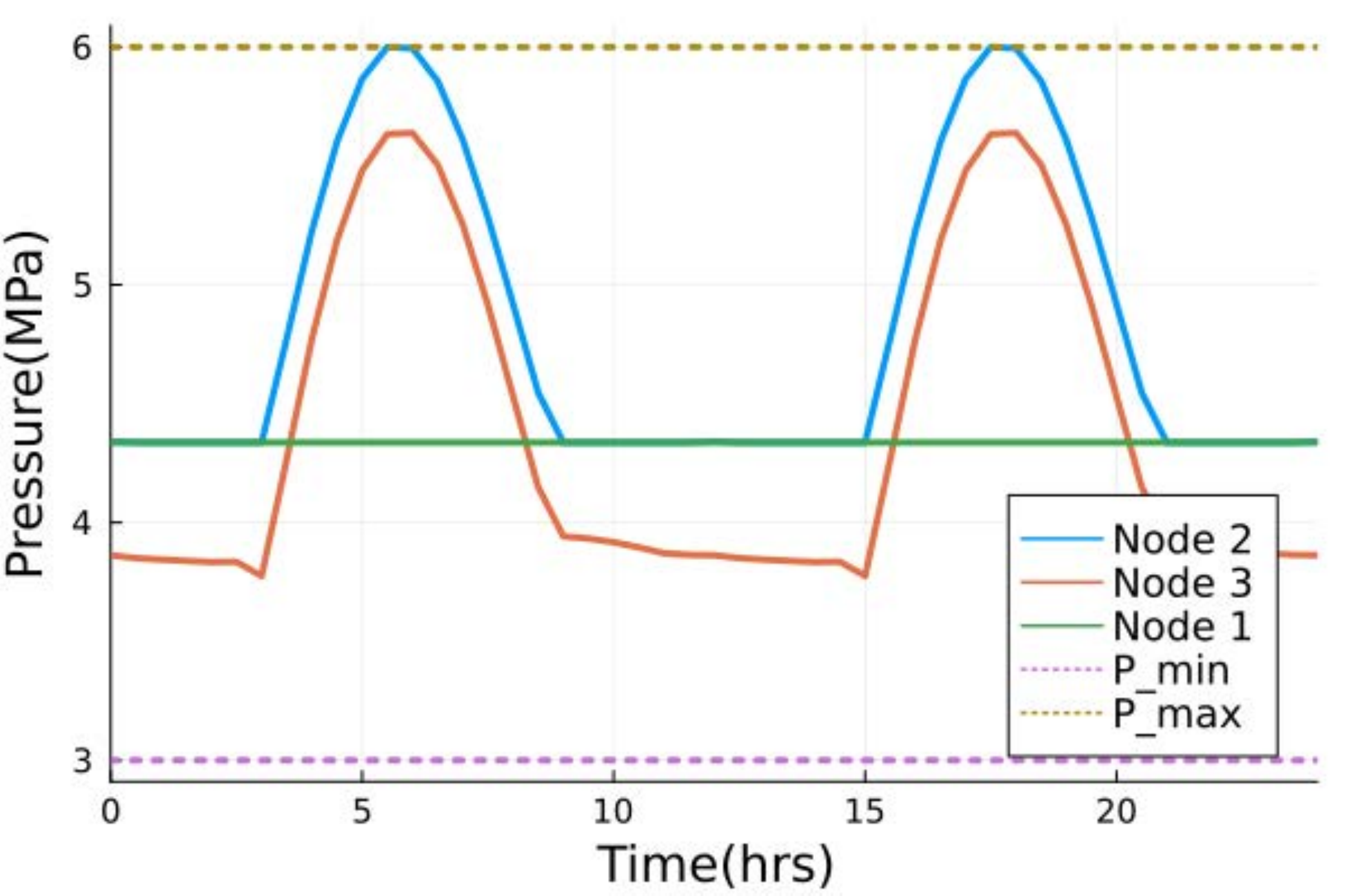}
        \caption{Nodal Pressure}
        \label{Node Pressure}
    \end{subfigure}
    \caption{Optimal pressure, density, and concentration solutions for the single pipe case study.}
    \label{single pipe plots}
    \vspace{-2ex}
\end{figure}

Consider a single pipe of length $L=30$ km, diameter $D=0.9144$ m, and friction factor $\lambda=0.01$ connected with a compressor near the injection node (1) as shown in Fig. \ref{sinlge-pipe}. We choose a discretization length of $\Delta L = 10$ km to formulate the lumped element reduced model. The time-varying injection concentration $\eta^s$ is given as
\begin{equation}\label{eq:node_injection_conc}
    \eta^s = \eta^s_0 + \delta \sin(2\pi\nu t/T),
\end{equation}
with $\eta^s_0 = 0.1$, $\delta = 0.05$, and $\nu = 2$. The function in equation \eqref{eq:node_injection_conc} is shown as the node 1 concentration in Fig. \ref{Node concentration}. Pressure at the slack injection node 1 is fixed at $p^s = 4.337$ MPa. The energy limit at the withdrawal node (3) is set at $g^{E,max} = 8000$ MJ/s, and we bound the flow through the compressor at $f^{c,\max} = 150$ kg/s. The pressure limits at the nodes are set at $p^{min} = 3.0$ MPa and $p^{max} = 6.0$ MPa.

The results from the optimization are plotted in Fig. \ref{single pipe plots} with profiles for nodal density, concentration, and pressure, and the boost ratio of the compressor.  The nodal values are of interest, because these are observed by operators of the pipeline and by system users at custody transfer points. We first observe that the profiles are periodic with two cycles ($\nu = 2.0$) in the 24 hour time period, and therefore discuss the trends in the first cycle ($t=0$ to $t=12$ hours) only. In Fig. \ref{Node concentration}, the \htwo concentrations at the three main nodes are plotted with time. The concentrations at node 1 (green) and 2 overlap because they are only connected by a compressor, whereas the concentration at the withdrawal node 3 (orange) is shifted to the right or lags behind because of the advective transport effect over the long pipe. The nodal density $\rho$ plotted in Fig. \ref{Node density} is the sum of partial densities of \htwo and NG. The density at node 1 (green) follows a smooth sinusoid similar to the injection node concentration because the pressure at node 1 (slack) is fixed and the two quantities are related by equation \eqref{eq:EOS2}. Nodal density depends on both the concentration and pressure at the node, which we observe in Figures \ref{Node concentration}, \ref{Node density}, and \ref{Node Pressure} where there is a sharp increase in compressor pressure outlet (node 2) around $t=3.5$ hours, which increases until it reaches the upper bound $p=6.0$ MPa at $t=5.5$ hours and then begins to decrease. In contrast to that, the densities at compressor outlet and withdrawal node vary slowly, and reached their maxima at $t=7.5$ hours (see Fig. \ref{Node density}). Observe also that the density at the withdrawal node (node 3) has a smaller local maximum near $t=9.5$ hours that results because the minimum in nodal concentration at the withdrawal node occurs at that time. The pressure drop across the pipe  seen in Fig. \ref{Node Pressure} appears higher when the node concentrations are high, such as at $t=2.5$ hours or when injection and withdrawal flows are high, as at $t=10.5$ hours.

\begin{figure}[t!]
    \centering
    \begin{subfigure}[b]{1.0\linewidth}
        \centering
        \includegraphics[height=40mm,width=0.9\linewidth]{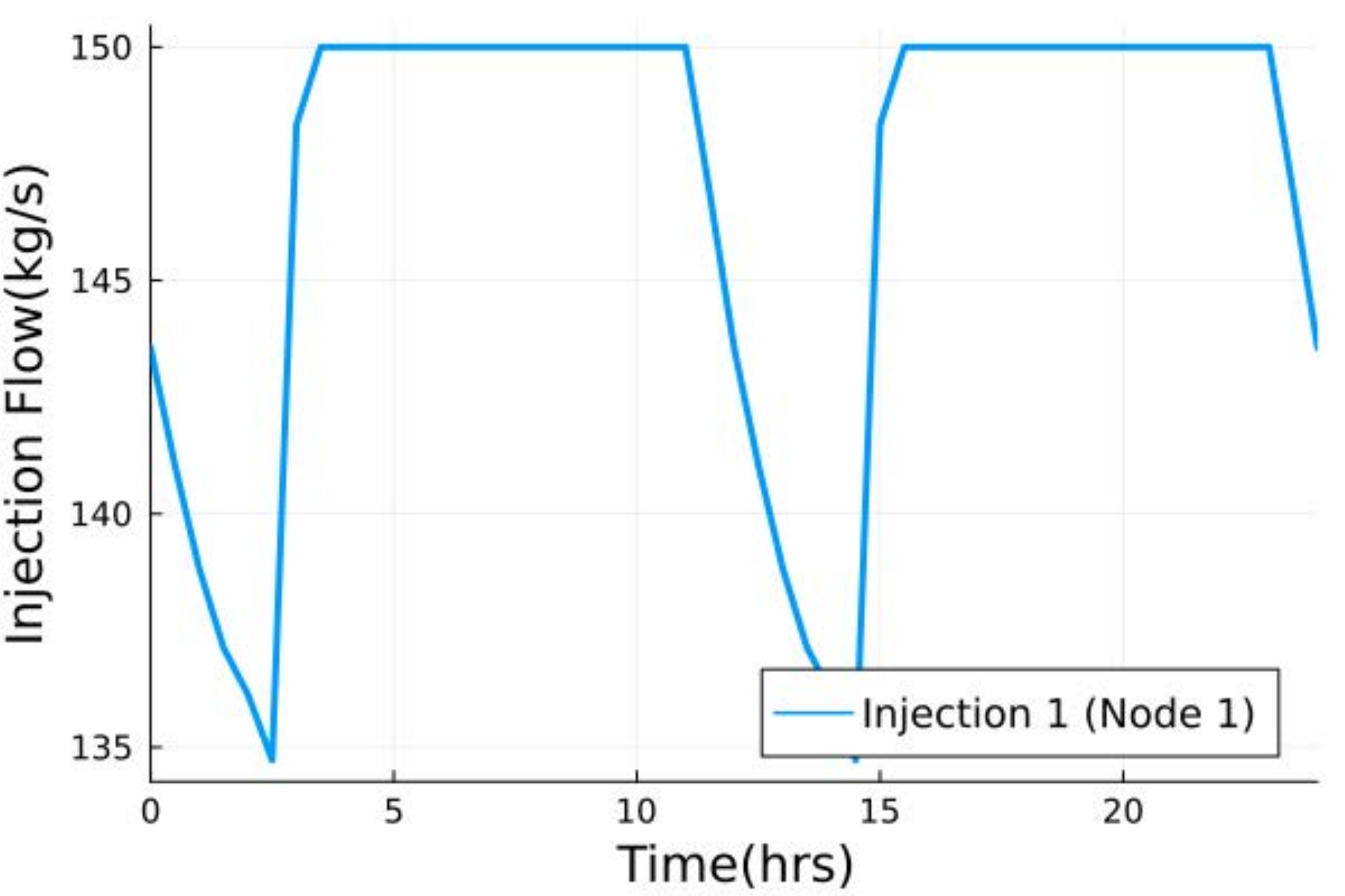}
        \caption{Injection flow}
        \label{Injection flow}
    \end{subfigure}
    \hfill
    \begin{subfigure}{1.0\linewidth}
        \centering
        \includegraphics[height=40mm,width=0.9\linewidth]{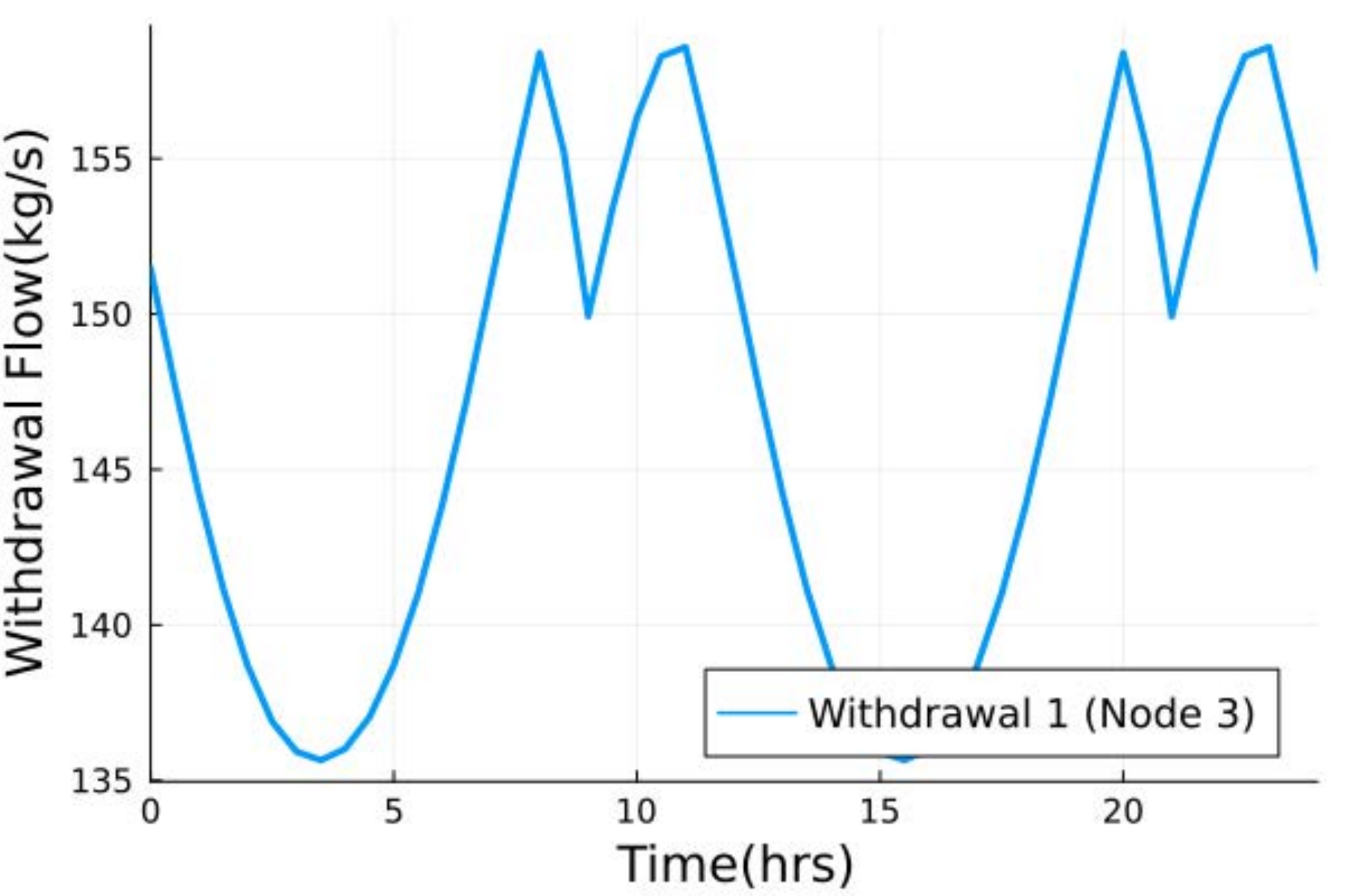}
        \caption{Withdrawal flow}
        \label{Withdrawal flow}
    \end{subfigure}
    \hfill
    \begin{subfigure}{1.0\linewidth}
        \centering
        \includegraphics[height=40mm,width=0.9\linewidth]{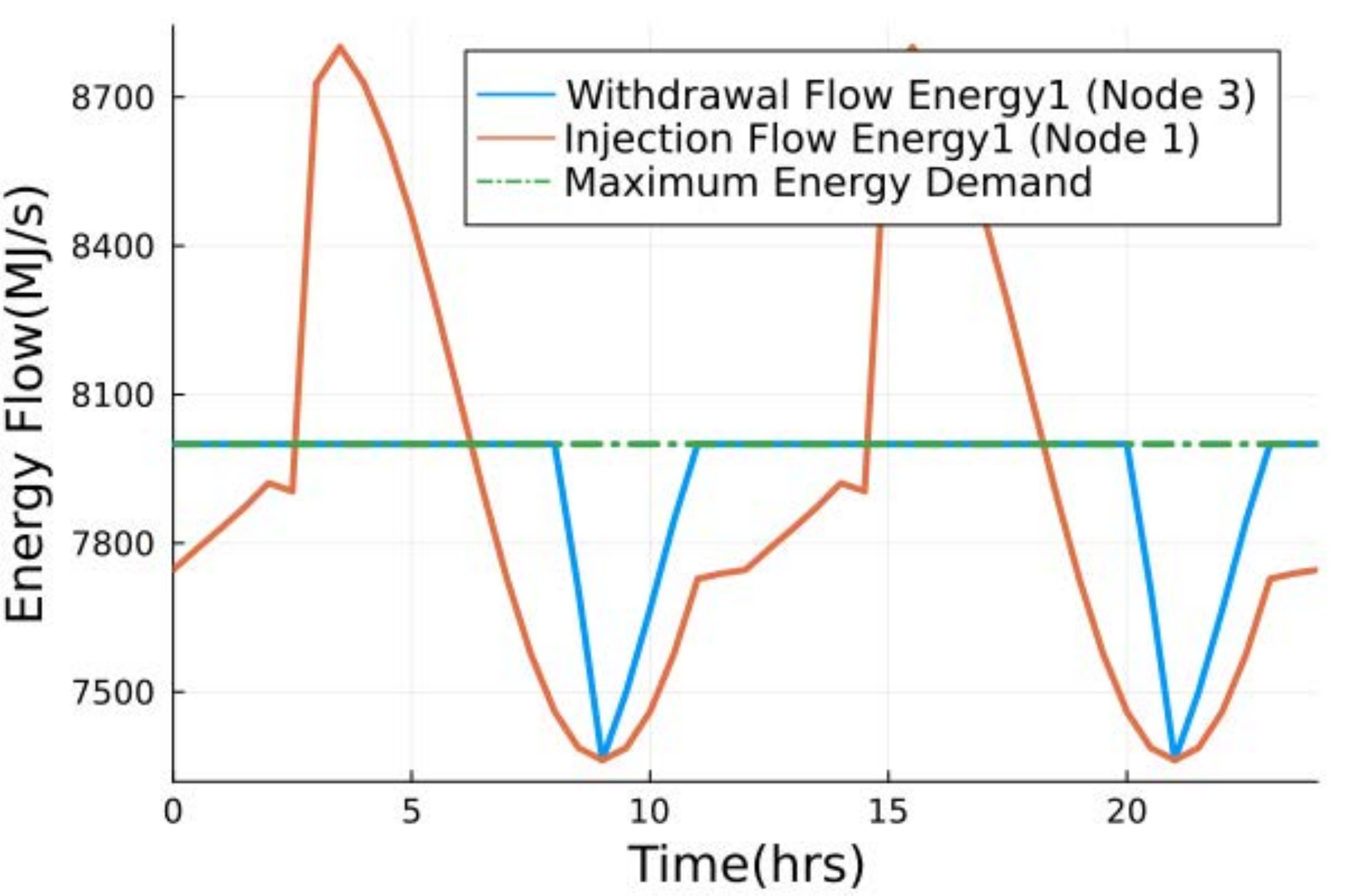}
        \caption{Energy Flow}
        \label{Withdrawal Energy}
    \end{subfigure}
    \caption{Optimal flow solutions for the single pipe case study.}
    \label{single pipe plots II}
    \vspace{-2ex}
\end{figure}

Fig. \ref{single pipe plots II} shows the variation in injection and withdrawal flow along with the energy delivered at the withdrawal node. Because the optimization solver seeks to maximize the energy delivered at the withdrawal node (see equation \eqref{eq:total_objective}), the amount of energy delivered is equal to its upper bound limit at 8000 MJ/s, except for a short duration between $t=9.0$ to $t=10.5$ hours. This is because the \htwo injection concentration in Fig. \ref{Node concentration} is so low that the energy content of the gas cannot meet the desired value of 8000 MJ/s without violating the compressor flow bound of 150 kg/s (see Fig. \ref{Withdrawal Energy}). Notice that the injection flow and energy flow both increase sharply at around $t=2.5$ hours, but the energy injection begins to decrease subsequently because of the decrease in injection concentration.  This is consistent with the pressure and node concentration profiles in Fig. \ref{Node Pressure} and Fig. \ref{Node concentration}, respectively. Similarly, the withdrawal flow in Fig. \ref{Withdrawal flow} is lower when nodal \htwo concentration is high and increases when the nodal concentrations begins to decrease at $t=3.5$ hours, until it is unable to match the maximum energy demand at $t=9.0$ hours, when it sharply decreases to minimize pressure drop, resulting in the compressor activity shown in Fig. \ref{Compressor Ratio}. The withdrawal flow subsequently begins to slowly increase again to maximize the energy supply when the withdrawal node concentration begins increasing at $t=10.0$ hours. Another important observation is that the maximum energy demand at the withdrawal node is met more regularly than the amount of time the energy flowing into node 1 is above the desired delivery rate at node 3.  

\subsection{8-node Network Case Study} \label{sec:eightnode}

\begin{figure}[t!]
    \centering
    \includegraphics[scale = 0.52]{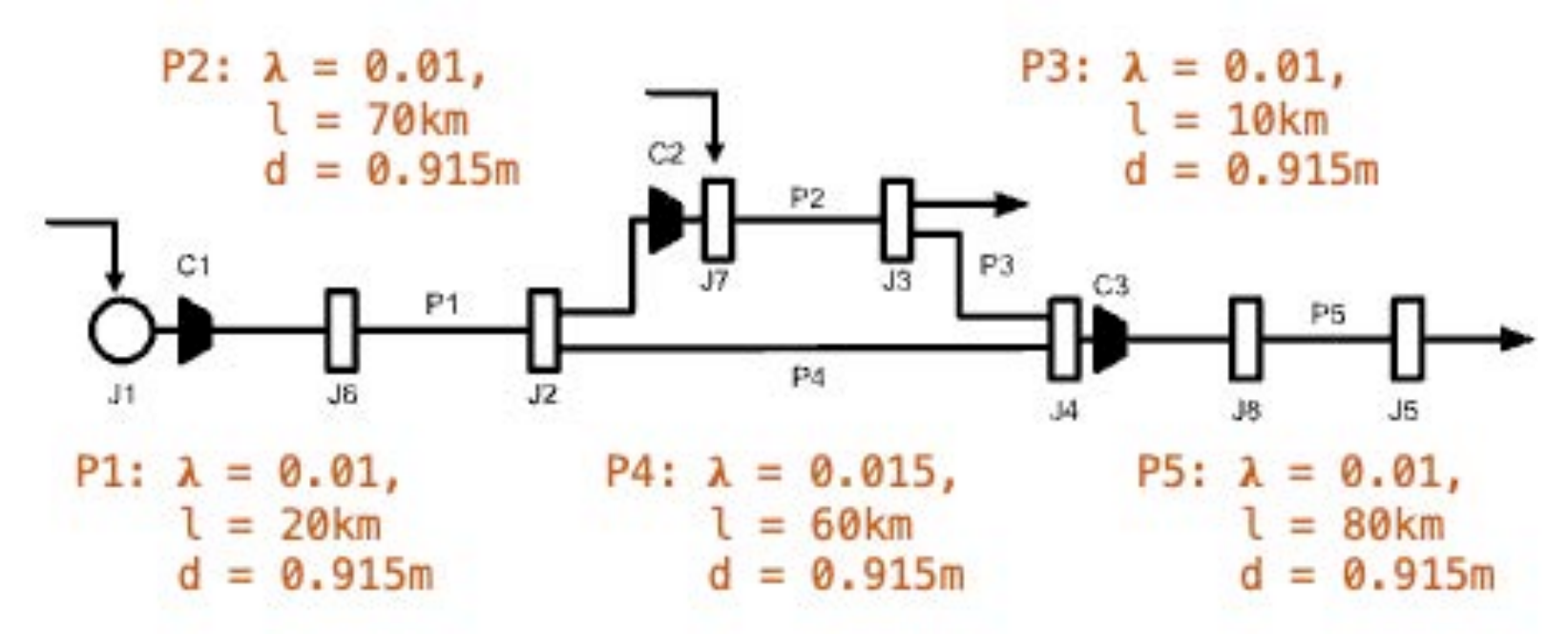}
    \caption{8-node test pipeline network and parameters.}
    \label{fig:8node}
    \vspace{-2ex}
\end{figure}

Next, we optimize an 8-node gas network with a loop, as shown in Figure \ref{fig:8node}.  The network has one slack node (J1), two injections (J1 and J7), and two withdrawals (J3 and J5) along with three compressors. The pressure and injection concentration at slack node J1 are fixed at $p=5.0$ MPa and $\eta^s = 0.08$. Meanwhile, the injection concentration at node J7 is varied as a sinusoid of the form in equation \eqref{eq:node_injection_conc}, similar to the single-pipe case, with $\eta^s_0 = 0.05$, $\delta = 0.01$, and $\nu = 1$. We use a coarser time discretization of $\Delta t = 1.0$ hour in this example, which results in 7536 variables and 7392 equality constraints. The compressor flow bounds are set at $f_1^{c,\max}=275$ kg/s, $f_2^{c,\max}=260$ kg/s and $f_3^{c,\max}=140$ kg/s. Pressure bounds and maximum energy demand are the same as in the first example.

\begin{figure}[th!]
    \centering
    \begin{subfigure}[b]{1.0\linewidth}
        \centering
        \includegraphics[height=40mm,width=0.9\linewidth]{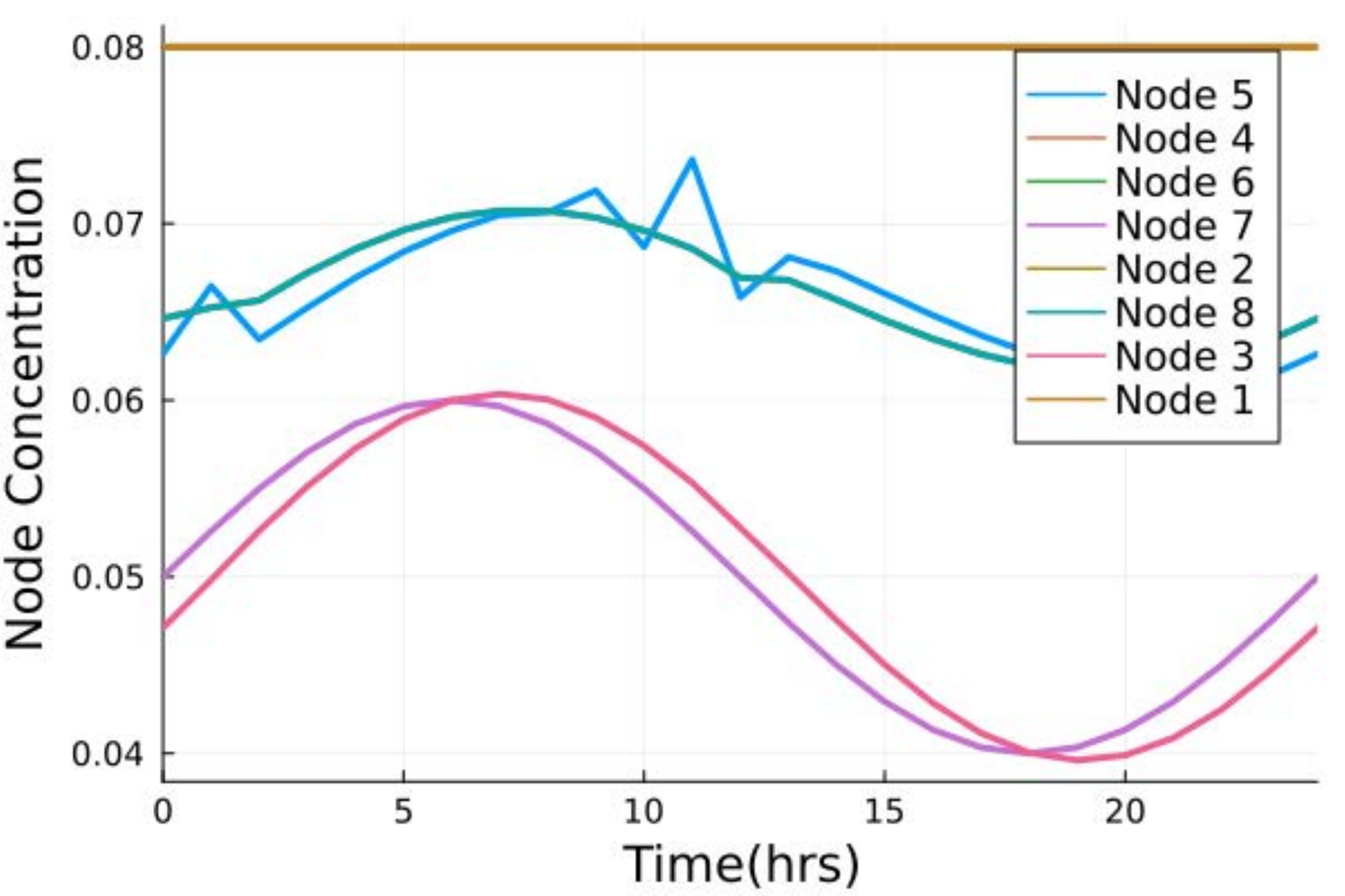}
        \caption{Nodal concentration}
        \label{8 node Node concentration}
    \end{subfigure}
    \hfill
    \begin{subfigure}[b]{1.0\linewidth}
        \centering
        \includegraphics[height=40mm,width=0.9\linewidth]{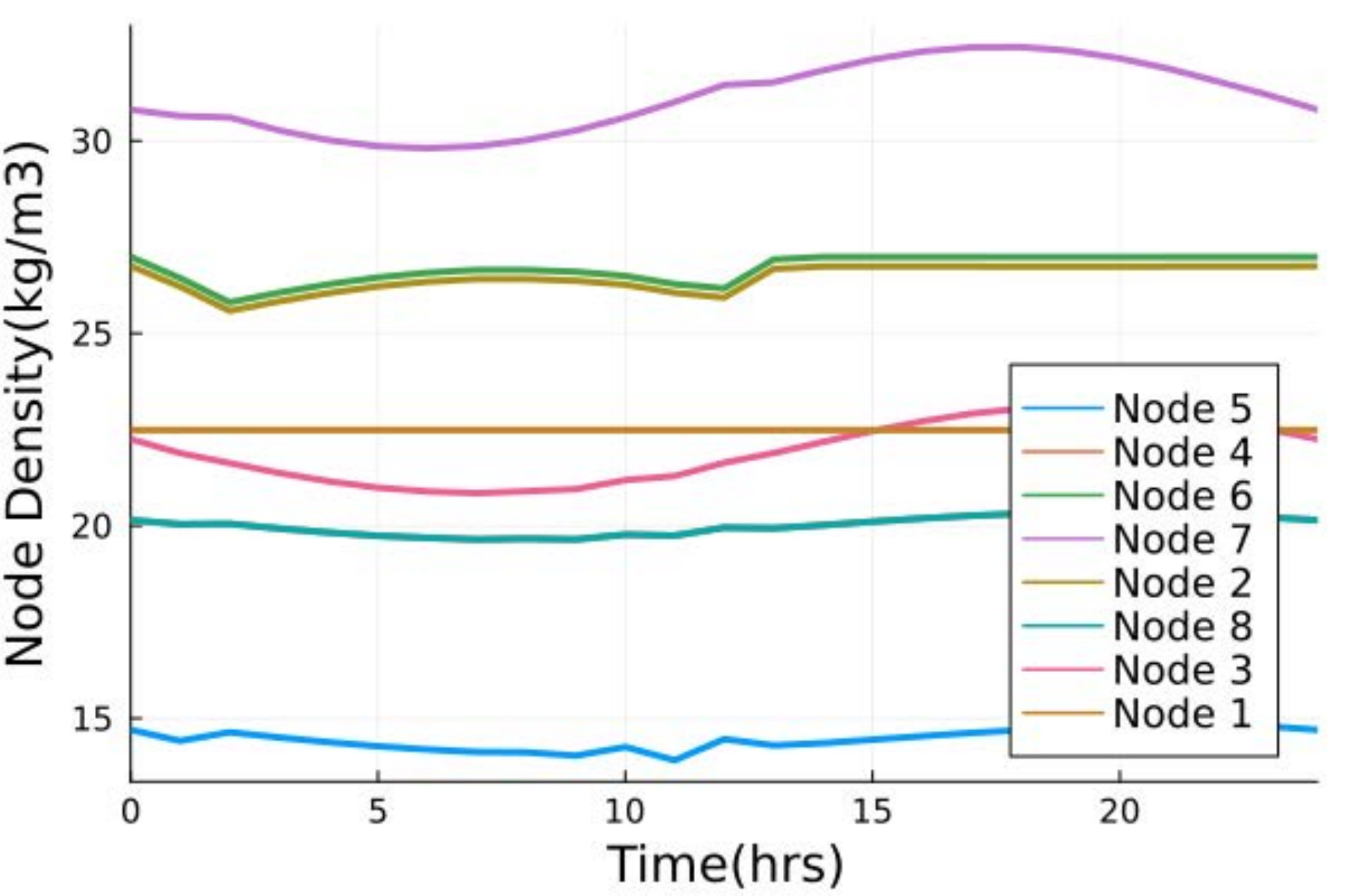}
        \caption{Nodal density}
        \label{8 node Node density}
    \end{subfigure}
    \hfill
    \begin{subfigure}[b]{1.0\linewidth}
        \centering
        \includegraphics[height=40mm,width=0.9\linewidth]{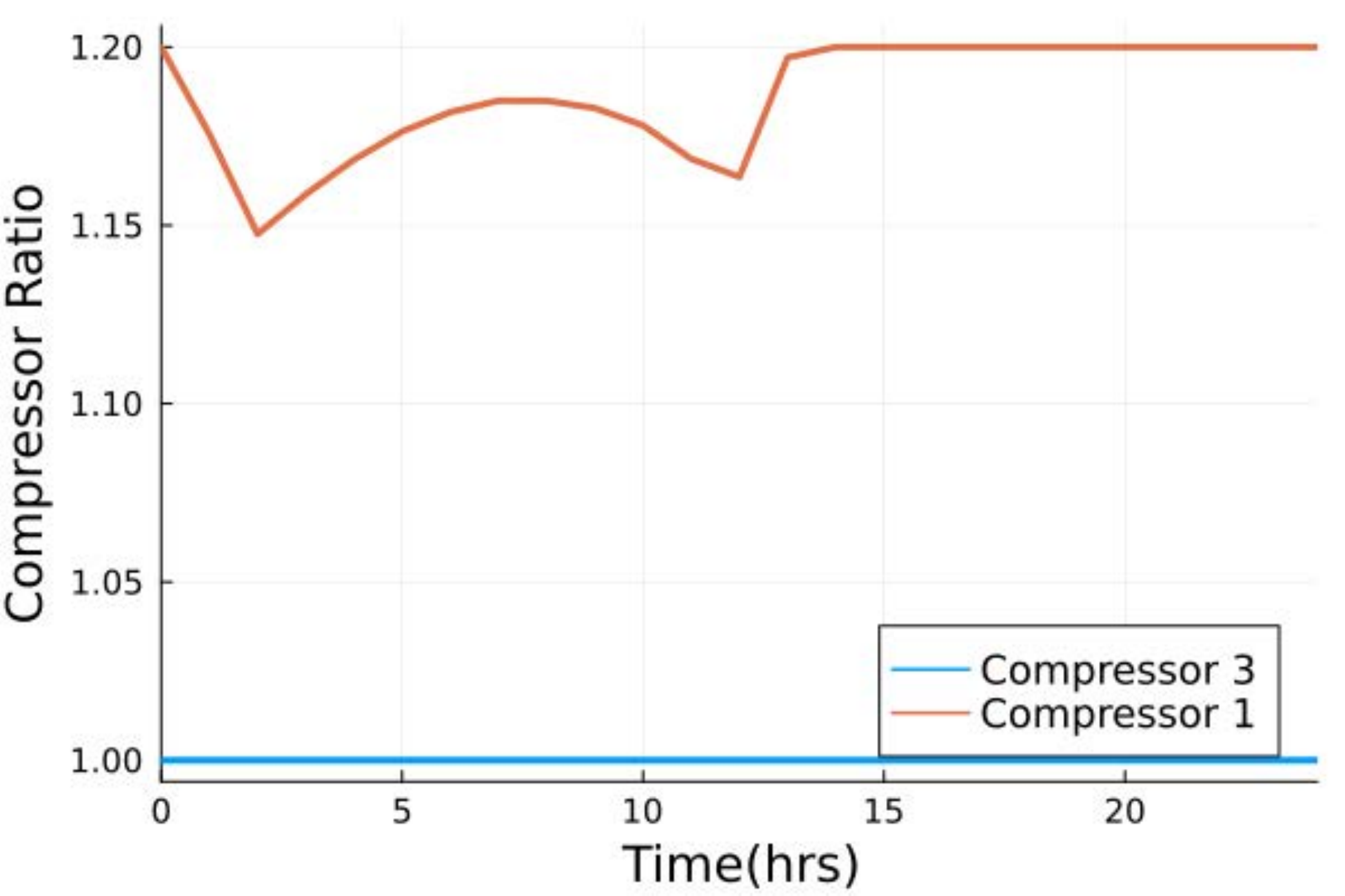}
        \caption{Compressor Ratio}
        \label{8 node Compressor Ratio}
    \end{subfigure}
    \hfill
    \begin{subfigure}[b]{1.0\linewidth}
        \centering
        \includegraphics[height=40mm,width=0.9\linewidth]{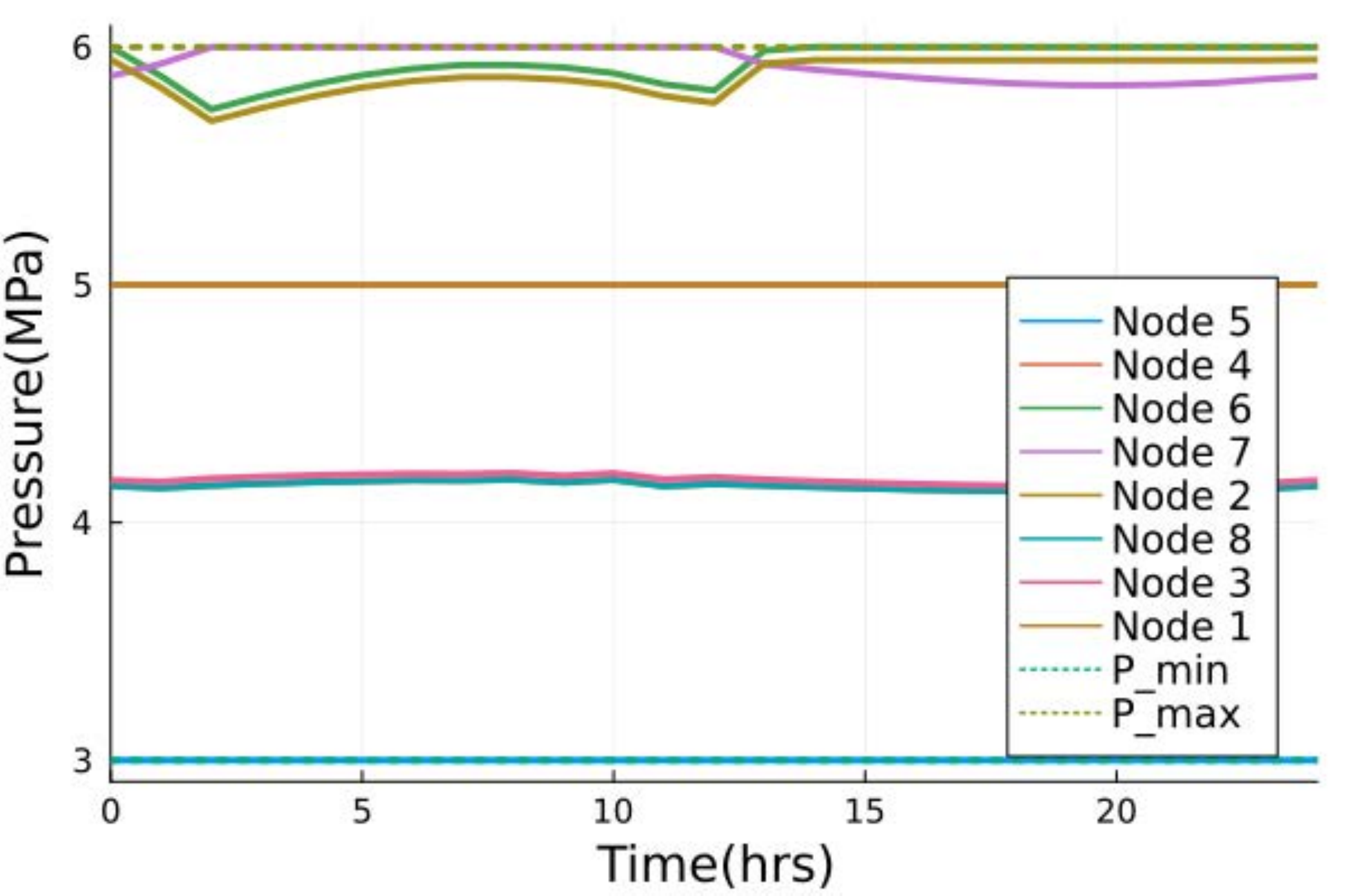}
        \caption{Nodal Pressure}
        \label{8 node Node Pressure}
    \end{subfigure}
    \caption{Optimal pressure, density, and concentration solutions for the 8-node test network.}
    \label{fig:8nodeplots}
    \vspace{-1ex}
\end{figure}

\begin{figure}[t!]
    \centering
    \begin{subfigure}[b]{1.0\linewidth}
        \centering
        \includegraphics[height=40mm,width=0.9\linewidth]{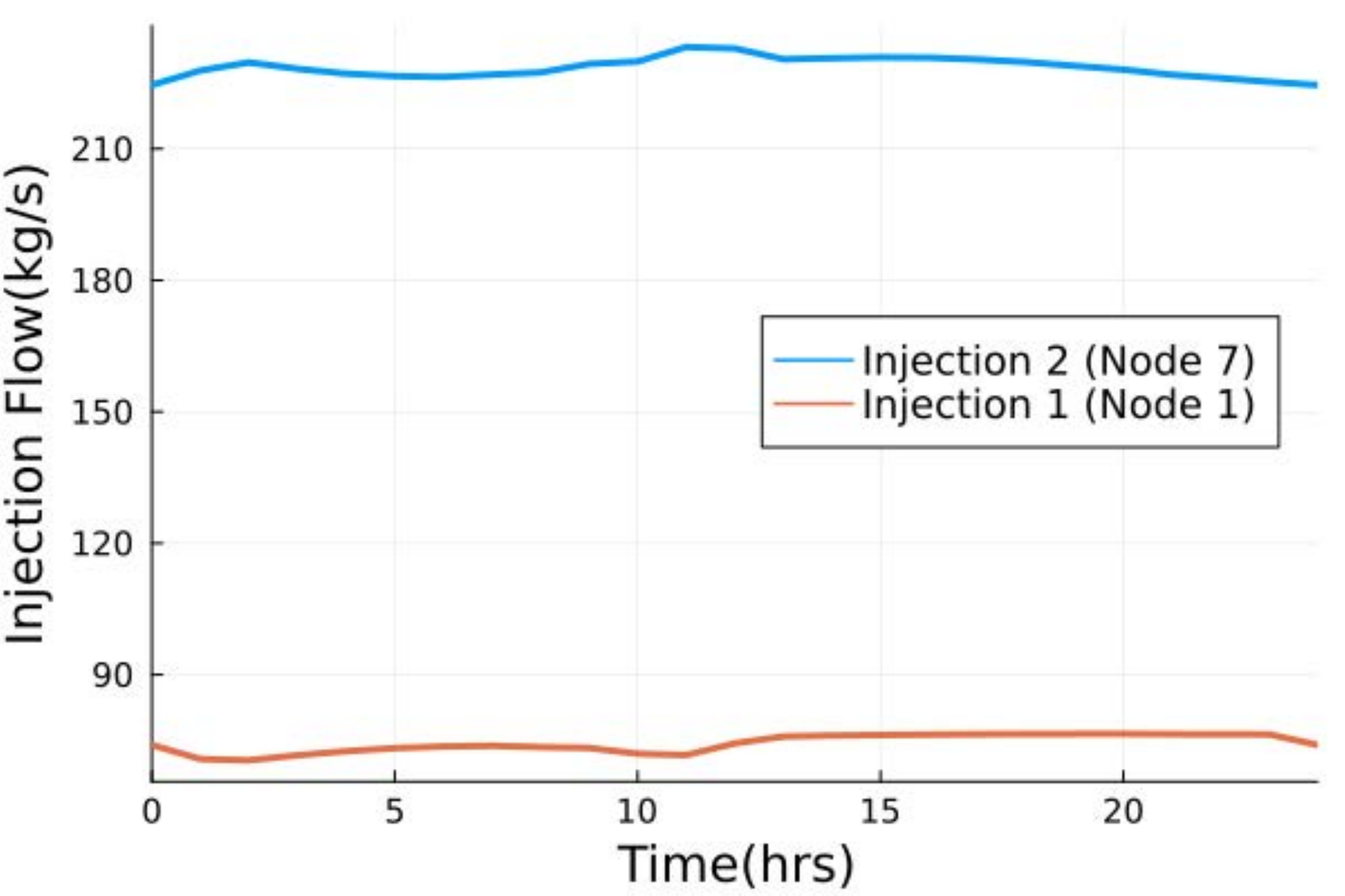}
        \caption{Injection flow}
        \label{8 node Injection flow}
    \end{subfigure}
    \hfill
    \begin{subfigure}{1.0\linewidth}
        \centering
        \includegraphics[height=40mm,width=0.9\linewidth]{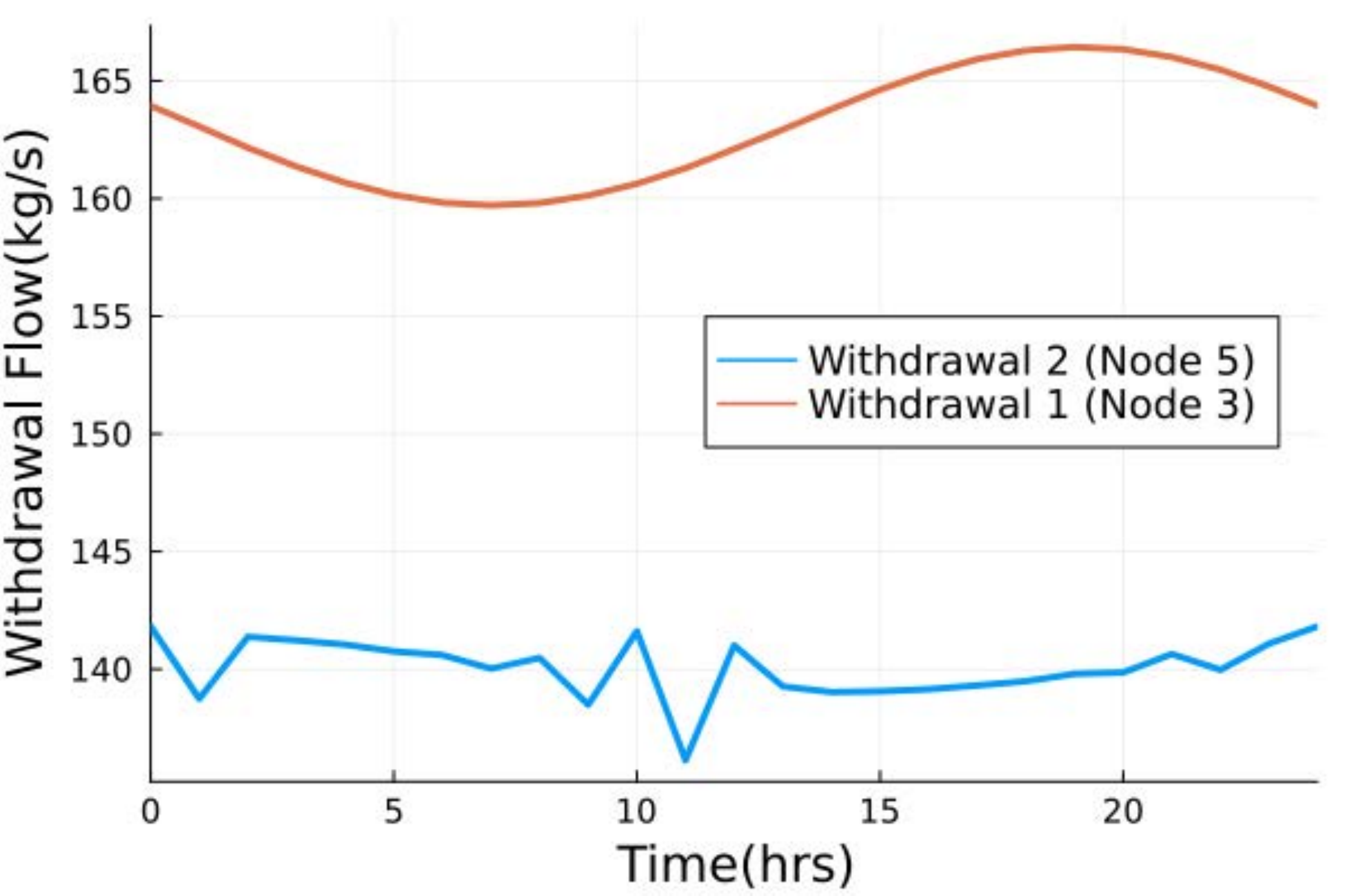}
        \caption{Withdrawal flow}
        \label{8 node Withdrawal flow}
    \end{subfigure}
    \hfill
    \begin{subfigure}{1.0\linewidth}
        \centering
        \includegraphics[height=40mm,width=0.9\linewidth]{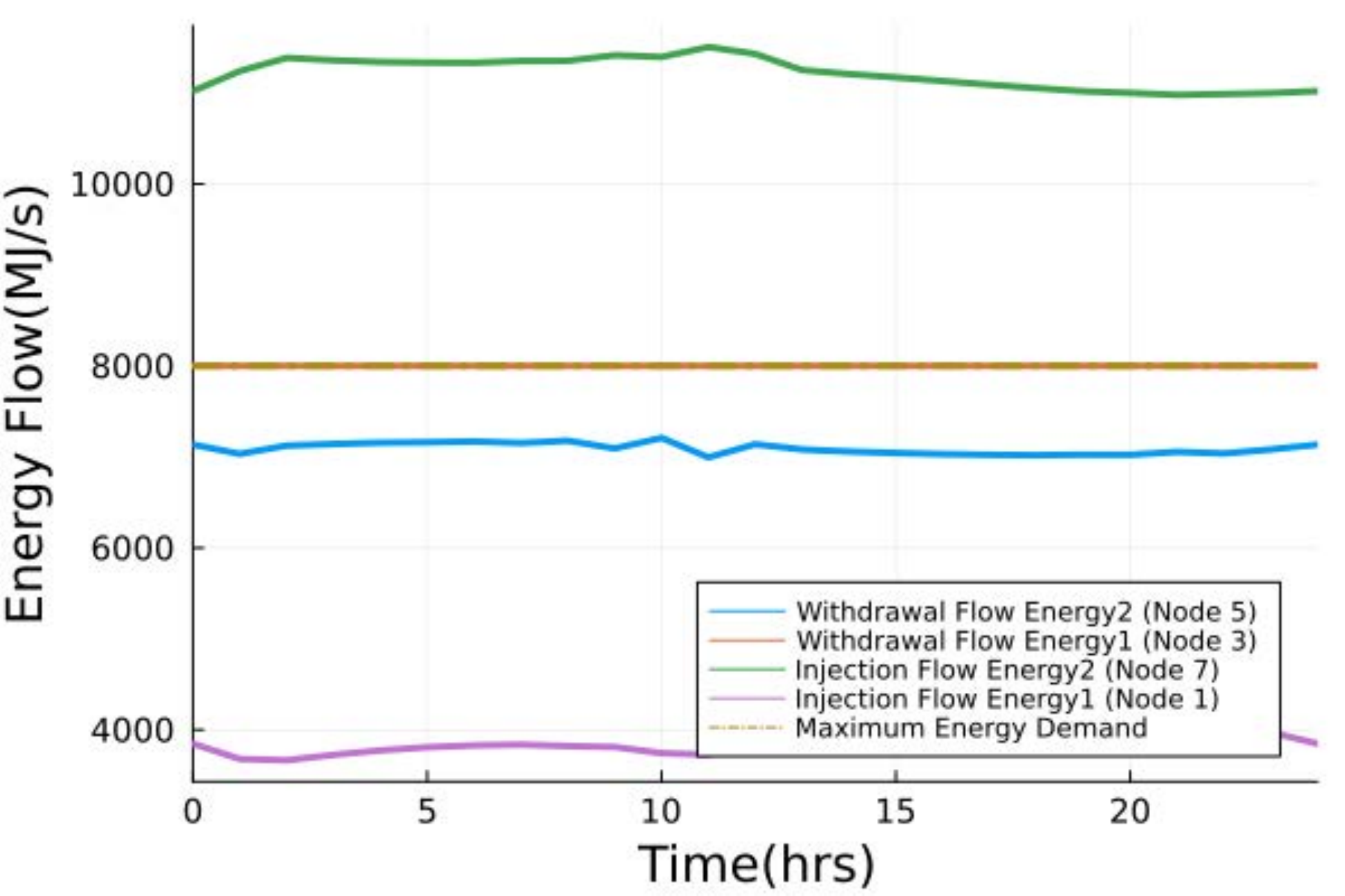}
        \caption{Energy Flow}
        \label{8 node Withdrawal Energy}
    \end{subfigure}
    \caption{Optimal flow solutions for the 8-node test network.}
	\vspace{2ex}
    \label{fig:8nodeplots2}
\end{figure}

The optimal solutions for the 8-node case study are shown in Fig. \ref{fig:8nodeplots} and Fig. \ref{fig:8nodeplots2}. Unlike for the single pipe case study, the profiles are less periodic and more non-uniform. A plausible explanation is that the nonlinearity of the nodal balance constraints in equations \eqref{eq:nodebalance} and the presence of two sources with different \htwo injection concentrations cause complex interactions with delays. Nonetheless, we do observe that there is no gas flow through compressor C2 in the optimal solution, and that compressor C3 is inactive in the optimal solution because the pressure at the withdrawal node J5 is at its lower bound. We also observe that the binding constraint in the network which is always active is the upper bound on the compressor flow C3, i.e., $f_3^{c,*} \equiv f_3^{c,\max}=140$ kg/s. The pressure at J6 and J2 follow a similar trend where the compressor C1 is used to increase the pressure of injected gas from J1 (see Fig. \ref{8 node Node Pressure} and Fig. \ref{8 node Compressor Ratio}). The \htwo concentration at node J4 (cyan) is determined by the nodal mixing of gas with fixed concentration from node J2 (yellow) and gas with varying concentration from node J3 (pink) as seen in Fig. \ref{8 node Node concentration}. In Fig. \ref{8 node Withdrawal flow}, the withdrawal flow at node J3 is sinusoidal, which balances the sinusoidal node concentration in Fig. \ref{8 node Node concentration} and results in constant withdrawal energy at 8000 MJ/s. 

\section{Conclusion and Future Work}    \label{sec:conclusions}

We formulate the challenging problem of model-predictive optimal control of hydrogen blending into natural gas pipeline networks subject to inequality constraints. Our numerical method employs a lumped parameter approximation to reduce the associated partial differential equations to a differential algebraic equation system that can be easily discretized and solved using nonlinear optimization solvers. We employ a circular time-discretization that is advantageous for time-periodic boundary conditions, parameters, and inequality constraint bound values.  The key advance in this study is the optimization of time-varying withdrawals from the network, where the running objective function is formulated in terms of energy delivery to consumers and mass injection by suppliers of natural gas or hydrogen.  This leads to a conceptually and computationally well-posed formulation, which enables numerical resolution of the complex behavior that arises from nonlinear interactions of inhomogeneous gas flows mixing throughout a pipeline network. Future studies will extend the setting to uncertain parameters and investigate economic questions by examining the dual solution.

\linespread{1}
\vspace{-0.5ex}
\bibliographystyle{unsrt} 
\bibliography{references}

\begin{thebibliography}{10}

\bibitem{chaczykowski2018gas}
Maciej Chaczykowski, Filip Sund, Pawe{\l} Zarodkiewicz, and Sigmund~Mongstad
  Hope.
\newblock Gas composition tracking in transient pipeline flow.
\newblock {\em J. Natural Gas Science \& Engineering}, 55:321--330, 2018.

\bibitem{baker2023transitions}
Luke~S. Baker, Saif~R. Kazi, and Anatoly Zlotnik.
\newblock Transitions from monotonicity to chaos in gas mixture dynamics in
  pipeline networks.
\newblock {\em PRX Energy}, 2:033008, Aug 2023.

\bibitem{tabkhi2008mathematical}
Firooz Tabkhi, Catherine Azzaro-Pantel, Luc Pibouleau, and Serge Domenech.
\newblock A mathematical framework for modelling and evaluating natural gas
  pipeline networks under hydrogen injection.
\newblock {\em International Journal of Hydrogen Energy}, 33(21):6222--6231,
  2008.

\bibitem{sodwatana2023optimization}
Mo~Sodwatana, Saif~R. Kazi, Kaarthik Sundar, and Anatoly Zlotnik.
\newblock Optimization of hydrogen blending in natural gas networks for carbon
  emissions reduction.
\newblock In {\em 2023 American Control Conference (ACC)}, pages 1229--1236.
  IEEE, 2023.

\bibitem{haeseldonckx2007use}
Dries Haeseldonckx and William D’haeseleer.
\newblock The use of the natural-gas pipeline infrastructure for hydrogen
  transport in a changing market structure.
\newblock {\em Internatl. J. Hydrogen Energy}, 32(10-11):1381--1386, 2007.

\bibitem{zlotnik2023effects}
Anatoly Zlotnik, Saif~R. Kazi, Kaarthik Sundar, Vitaliy Gyrya, Luke Baker,
  Mo~Sodwatana, and Yan Brodskyi.
\newblock Effects of hydrogen blending on natural gas pipeline transients,
  capacity, and economics.
\newblock In {\em PSIG Annual Meeting}, page 2312. PSIG, 2023.

\bibitem{guandalini2017dynamic}
Giulio Guandalini, Paolo Colbertaldo, et~al.
\newblock Dynamic modeling of natural gas quality within transport pipelines in
  presence of hydrogen injections.
\newblock {\em Applied Energy}, 185:1712--1723, 2017.

\bibitem{Kai_and_Saif}
Kai Liu, Saif~R. Kazi, Lorenz~T. Biegler, et~al.
\newblock Dynamic optimization for gas blending in pipeline networks with gas
  interchangeability control.
\newblock {\em AIChE Journal}, 66(5):e16908, 2020.

\bibitem{melaina2013blending}
Marc~W. Melaina, Olga Antonia, and Michael Penev.
\newblock Blending hydrogen into natural gas pipeline networks: A review of key
  issues.
\newblock {\em NREL Tech. Rep.}, NREL/TP-5600-51995, 2013.

\bibitem{miao2021long}
Bin Miao, Lorenzo Giordano, and Siew~Hwa Chan.
\newblock Long-distance renewable hydrogen transmission via cables and
  pipelines.
\newblock {\em International Journal of Hydrogen Energy}, 46(36):18699--18718,
  2021.

\bibitem{schuster2020centrifugal}
Sebastian Schuster, Hans~Josef Dohmen, and Dieter Brillert.
\newblock Challenges of compressing hydrogen for pipeline transportation with
  centrifugal compressors.
\newblock In {\em Proceedings of Global Power and Propulsion Society}, pages
  2504--4400, 2020.

\bibitem{zhang2022modelling}
Zihang Zhang, Isam Saedi, Sleiman Mhanna, Kai Wu, and Pierluigi Mancarella.
\newblock Modelling of gas network transient flows with multiple hydrogen
  injections and gas composition tracking.
\newblock {\em International Journal of Hydrogen Energy}, 47(4):2220--2233,
  2022.

\bibitem{witkowski2018analysis}
Andrzej Witkowski, Andrzej Rusin, Miros{\l}aw Majkut, and Katarzyna Stolecka.
\newblock Analysis of compression and transport of the methane/hydrogen mixture
  in existing natural gas pipelines.
\newblock {\em International Journal of Pressure Vessels and Piping},
  166:24--34, 2018.

\bibitem{rachford2000optimizing}
Henry~H. Rachford and Richard~G. Carter.
\newblock Optimizing pipeline control in transient gas flow.
\newblock In {\em PSIG Annual Meeting}, page 0004, 2000.

\bibitem{hari2021operation}
Sai Krishna~Kanth Hari, Kaarthik Sundar, Shriram Srinivasan, Anatoly Zlotnik,
  and Russell Bent.
\newblock Operation of natural gas pipeline networks with storage under
  transient flow conditions.
\newblock {\em IEEE Transactions on Control Systems Technology},
  30(2):667--679, 2021.

\bibitem{osiadacz1984simulation}
Andrzej Osiadacz.
\newblock Simulation of transient gas flows in networks.
\newblock {\em International J. Numerical Methods in Fluids}, 4(1):13--24,
  1984.

\bibitem{rudkevich2019evaluating}
Aleksandr Rudkevich, Anatoly Zlotnik, Xindi Li, Pablo Ruiz, et~al.
\newblock Evaluating benefits of rolling horizon model predictive control for
  intraday scheduling of a natural gas pipeline market.
\newblock In {\em 52nd Hawaii International Conference on System Sciences},
  pages 3627--3636, 2019.

\bibitem{herty2010new}
Michael Herty, Jan Mohring, and Veronika Sachers.
\newblock A new model for gas flow in pipe networks.
\newblock {\em Mathematical Methods in the Applied Sciences}, 33(7):845--855,
  2010.

\bibitem{grundel2014model}
Sara Grundel, Lennart Jansen, Nils Hornung, et~al.
\newblock Model order reduction of differential algebraic equations arising
  from the simulation of gas transport networks.
\newblock In {\em Progress in Differential-Algebraic Equations: Deskriptor
  2013}, pages 183--205. Springer, 2014.

\bibitem{Sundar2018}
Kaarthik Sundar and Anatoly Zlotnik.
\newblock State and parameter estimation for natural gas pipeline networks
  using transient state data.
\newblock {\em IEEE Transactions on Control Systems Technology},
  27(5):2110--2124, 2018.

\bibitem{menon05}
E.~S. Menon.
\newblock {\em Gas Pipeline Hydraulics}.
\newblock CRC Press, 2005.

\bibitem{knitro}
Richard~H. Byrd, Jorge Nocedal, et~al.
\newblock {\em Knitro: An Integrated Package for Nonlinear Optimization}, pages
  35--59.
\newblock Springer, Boston, 2006.

\bibitem{jump}
Iain Dunning, Joey Huchette, and Miles Lubin.
\newblock Jump: A modeling language for mathematical optimization.
\newblock {\em SIAM Review}, 59(2):295--320, 2017.

\end{thebibliography}






\end{document}